\documentclass[11pt]{article}
\usepackage[dvips]{graphics}
\usepackage{amssymb,latexsym}

\newcommand{\R}{{\bf R}}

\newcommand{\ga}{{\gamma}}
\newcommand{\Ga}{{\Gamma}}

\newcommand{\la}{{\lambda}}

\newcommand{\Cc}{{\cal C}}
\newcommand{\Dd}{{\cal D}}

\newcommand{\Jj}{{\cal J}}

\newcommand{\Uu}{{\cal U}}

\newcommand{\proof}[1]{\noindent{\bf Proof#1:\  }}

\newcommand{\QED}{\hfill$\Box$\medskip}

\newcommand{\CP}{{\rm CP^2}}
\newcommand{\CPb}{{\rm \widetilde{\CP}}}
\newcommand{\cHZ}{{c_{HZ}}}
\newtheorem{theorem}{Theorem}[section]

\newtheorem{definition}[theorem]{Definition}

\newtheorem{lemma}[theorem]{Lemma}

\newtheorem{prop}[theorem]{Proposition}
\newtheorem{proposition}[theorem]{Proposition}
 
\newtheorem{makingmainthmgeod}{Theorem \ref{theorem: mainthmgeod}}

\newtheorem{makingchzprod}{Theorem \ref{theorem: chzprod}}

\title {\bf Hofer-Zehnder capacity and length minimizing paths in the Hofer norm}
\author{\textsc{ JENNIFER SLIMOWITZ}}
\date{\textsc{\today}}

\begin{document}

\maketitle 

\abstract{ We use  the criteria of 
Lalonde and McDuff to determine a new class of examples of length 
minimizing paths in the group $Ham(M)$.  
For a compact symplectic manifold $M$ of dimension two or four, we show that a path 
in $Ham (M)$, generated by an autonomous Hamiltonian and starting at the identity, which induces
no non-constant closed trajectories of points in $M$, is length minimizing
among all homotopic paths. The major step in the proof involves determining  
an upper bound for the Hofer-Zehnder capacity for symplectic manifolds of the type
$(M \times D(a))$ where $M$ is compact and has dimension two or four. 
In the appendix, we give an alternate proof
of Polterovich's  result that rotation in $\CP$ and in the blow-up 
of $\CP$ at one point is a length minimizing  path with respect 
to the Hofer norm.  Here we use the Gromov capacity and describe the
necessary ball embeddings. }

\section{Background and Main Theorems}

In this paper, we show that certain naturally occurring paths in the
group $Ham(M)$ are length minimizing  paths with respect to the
Hofer norm. 
A length minimizing  path $\phi_t$ for $0 \leq t \leq 1$ in $Ham(M)$
is a path which is an absolute minimum of the length functional among all paths 
from $\phi_0$ to $\phi_1$.   The search for length minimizing  paths is the
logical extension of the work done on general geodesics, that is those which
minimize length locally, by Bialy-Polterovich in \cite{BP}, Ustilovsky in \cite{U},
and Lalonde-McDuff in \cite{LM2}.   
 
Let $(M, \omega)$ be any symplectic manifold and $H_t : M \rightarrow {\bf R}$
for $0 \leq t \leq 1$ be a compactly supported time-dependent 
Hamiltonian function on $M$.  The length $L(H)$ of $H$ is defined to be
$$L(H) = \int_0 ^1 \max_{x \in M} H_t(x) - \min_{x \in M} H_t(x) dt.$$
 The time-dependent Hamiltonian vector field $X^H$ induced by $H$ is the unique solution

\smallskip 
\noindent\rule{1in}{.2mm}

\noindent{\footnotesize \textit{Mathematics Subject Classification:} 53C15, 58F05, 58D05, 58B20 }

\nopagebreak
\noindent{\footnotesize\textit{Key words and phrases:} symplectic geometry, Hamiltonian diffeomorphisms, Hofer norm,
 Hofer-Zehnder capacity}

\noindent to the equation
$$i(X^H) \omega = - dH,$$
and its time $t$ flow is denoted $\phi ^H _t$.  
The group $Ham^c(M)$ is the set of compactly supported time-one Hamiltonian maps on $M$:
$$Ham^c(M) = \{ \phi ^H _1 \; \vline \; H_t : M \rightarrow {\bf R} \}. $$
When working on a compact manifold, we drop the superscript $c$ and write only
$Ham(M)$ since all diffeomorphisms have compact support.

Now, to any path $\phi _t$ for $0 \leq t \leq 1$ in $Ham^c(M)$, 
we can associate its generating time-dependent Hamiltonian $H_t$ satisfying
$$\frac{d}{dt}\vline_{t=t_0} \phi_t = X^H_{t_0}(\phi ^H _{t_0}).$$ The
length of the path $\phi _t$ is defined as the length $L(H)$ of its generating
Hamiltonian.   The Hofer norm $\parallel \phi \parallel$ of 
$\phi \in Ham^c(M)$ is
the infimum of the lengths of all of the paths from the identity to $\phi$;
hence, a globally length minimizing path from the identity to $\phi$
determines $\parallel \phi \parallel$.  Although the Hofer
norm is simply defined, it is difficult to calculate.  One case in which
it might be easiest to calculate $\parallel \phi \parallel$ is when there
is a natural path from the identity to $\phi$, e.g. a path induced by 
a circle action such as a rotation.

Lalonde and McDuff provide an important example of a globally
length minimizing path when they show that rotation through $\pi$ radians
on $S^2$ is  length minimizing in $Ham(S^2)$ \cite{LM2}.  
This leads us to ask whether rotation of $\CP$ through $\pi$ radians is length 
minimizing in $Ham(\CP)$.  In fact, by following the procedure outlined by Lalonde and McDuff
in \cite{LM2} using  quasi-cylinders and capacities, we prove in the appendix that 
rotation  on $\CP$ and the blow-up of $\CP$ at one point is indeed a length minimizing  path.   
We work with Gromov capacity and show the necessary criteria are satisfied by constructing 
explicit embeddings of balls. These are  independent proofs 
of results that  Polterovich derives in \cite{P1} and \cite{P2}.

However, the power of Gromov capacity to detect length minimizing paths is limited,
and to obtain more general results we switch gears and work with the Lalonde-McDuff
criteria paired with the Hofer-Zehnder capacity instead of the Gromov capacity.  
Before stating results, we need  the follwing definition.  

A path $\phi_t \in Ham (M)$ which starts from the identity has 
{\bf no non-constant closed trajectory in time less than 1} if 
$$\phi_{t_0} (x_0) = x_0 \mbox{ for some } t_0 
\in (0, 1], x_0 \in M \Rightarrow \phi_t(x_0) = x_0 \; \forall t \in [0,1].$$
The next theorem is the main  result of this paper;  note that the results 
about rotation in $\CP$ and its
blow-up from the appendix can also be derived as an application of this theorem.  

\begin{makingmainthmgeod}
Let $(M, \omega)$ be a compact symplectic manifold of dimension two or four.
 Let $\phi^H_t$ for
$0 \leq t \leq 1$ be a path in $Ham (M)$  generated by an autonomous 
Hamiltonian 
$H: M \rightarrow {\bf R}$  such that $\phi^H_0$ is the identity
diffeomorphism and  $\phi^H_t$ has no non-constant closed trajectory in
time less than 1.
Then, the path $\phi^H_t$ for $0 \leq t \leq 1$ is length minimizing among all 
homotopic paths between the identity and $\phi^H_1$.
\end{makingmainthmgeod}

Theorem \ref{theorem: mainthmgeod} generalizes Hofer's parallel result for
 ${\bf R^{2n}}$.  His proof  that the flow of an autonomous
Hamiltonian in ${\R^{2n}}$  which admits no non-constant closed
trajectory in time less than 1 is a length minimizing path appears in Section 5.7
of \cite{HZ}.  In addition, Theorem \ref{theorem: mainthmgeod} is an
extension of  Lalonde and McDuff's Theorem \nolinebreak5.4 from \cite{LM2}.   There they
show  that the conclusion holds if $M$ has dimension two or if
$M$ is weakly exact.  Siburg has generalized Hofer's result in another way;
in \cite{Si} he extends the class of allowable Hamiltonians on ${\R^{2n}}$ to
include time dependent functions as well as autonomous ones.

By the classification paper of Karshon, we know exactly what the semi-free
Hamiltonian $S^1$ actions on a compact symplectic four manifold look like \cite{K}.  
Hence, if $M$ has dimension four and $H$ actually generates a loop, i.e. the path $\phi^H_t$ 
represents a circle action, we know up to  equivariant isomorphism 
the possible ways in which a $\phi^H_t$ that satisfies the hypotheses of Theorem 
\ref{theorem: mainthmgeod} rotates $M$.  

For the proof of Theorem
\ref{theorem: mainthmgeod}, we follow the criteria for length
minimizing  paths from \cite{LM2} using the  Hofer-Zehnder capacity.
Let $D(a)$ denote the open two-disk
equipped with a symplectic form $\sigma$
which satisfies $\int_{D} \sigma = a$.  In order to complete the proof,
we need to  show that the Hofer-Zehnder capacity $\cHZ$ satisfies
the capacity-area inequality on all manifolds of the form $M \times D(a)$,
equipped with the product symplectic form,
where $M$ is a symplectic manifold of dimension two
or four.  In \cite{HV}, Hofer and Viterbo  
have proven that $\cHZ$ satisfies this inequality for all $a >0$
if the manifold $M$ is weakly
exact.  This is a very restrictive condition which in particular excludes
the case $M = \CP$ or the blow-up of $\CP$.  Hence, in Section 4, we return 
to the original proof of Hofer and Viterbo in \cite{HV} and modify it using
the theory of J-holomorphic curves, proving:

\begin{makingchzprod}
 Let  $(M, \omega)$  be a compact symplectic manifold of dimension two or four.  Then,
$$\cHZ (M \times D(a), \omega \oplus \sigma) \leq a.$$
\end{makingchzprod}

\noindent {\bf Remark}
Theorems \ref{theorem: mainthmgeod} and \ref{theorem: chzprod}
as they are now stated have limited scope.  The restriction to manifolds of
dimension two or four is required in order to deal with multiply covered curves
on $M \times S^2$  at the end of Section 4.  However, recent advances in the
theory of $J-$holomorphic curves by Fukaya-Ono, Li-Tian,
Liu-Tian, McDuff, Ruan, and Siebert, following ideas of Konstevitch,  will most likely allow us to
generalize to other dimensions.  In particular, the methods that Liu and Tian  
use in \cite{LT} and McDuff's work in \cite{McD} that deal with stable virtual
moduli spaces of curves can probably be used to extend
Theorem \ref{theorem: chzprod} and   Theorem \ref{theorem: mainthmgeod}
to include manifolds of all dimensions.

\medskip

In  related work, Polterovich examines a rotation, similar to
the one considered in the appendix, on $\CP$ and on the monotone 
manifold $(\CPb, \tau_{1/\sqrt{3}})$, 
the blow-up of $\CP$ obtained by removing a ball of radius $\frac{1}{\sqrt{3}}$ centered at 
the point $[1:0:0]$.  He examines the path $\psi_t$ where
$$\psi_t [z_0: z_1: z_2] = [e^{2\pi it}z_0: z_1: z_2].$$ 
He shows that the loop formed by $\psi_t$ for $0 \leq t \leq 1$ 
is a length minimizing representative of its homotopy class in $Ham (\CP)$
in \cite{P1} and in $Ham (\CPb)$ in \cite{P2}.   Note that his results in \cite{P1}
and \cite{P2}  imply Theorems \ref{theorem: rotCP2} and \ref{theorem: rotCPtwoB} of this
paper;  however his proofs rely on  Gromov's K-area and a homomorphism combining
the symplectic action and the Maslov index. The proofs here using symplectic capacities
and quasi-cylinders illustrate the criteria described in \cite{LM2}. 

This paper is organized in the following way.  
The second section  describes the criteria for length minimizing
paths developed by Lalonde and McDuff in \cite{LM2}.
The third and fourth sections use  $J-$holomorphic curve theory to prove  
Theorem \ref{theorem:  mainthmgeod} and Theorem \ref{theorem: chzprod}.
The appendix of this paper gives in full detail the ball embeddings which show that
 specific rotations on $\CP$ and $\CPb$ are length minimizing . 

I thank my advisor Dusa McDuff for introducing me to the questions in
this paper, offering numerous suggestions, and reading several previous
drafts.  In addition, I thank Leonid Polterovich for his helpful comments
and insight and for pointing out several errors in a previous draft.   I am grateful to 
Francois Lalonde for  explaining key points and providing assistance with the bubbling
arugments.   This material is based upon work supported by the North
Atlantic Treaty Organization under a Grant awarded in 1998.

\section{Criteria for length minimizing paths}

We now briefly describe the theory that Lalonde and McDuff use to develop
their criteria for length minimizing  paths.  In \cite{LM2}, they first
derive a geometric way of detecting that $L(H_t) \leq L(K_t)$ for two Hamiltonians $H_t$ and
$K_t$ on $M$ .  Then, they determine sufficient 
conditions involving symplectic capacities for this geometric requirement
to be satisfied.    We also describe the Gromov 
capacity and the Hofer-Zenhder capacity, two symplectic capacities used in this paper.  

\subsection{Results of Lalonde and McDuff}
 To begin, we must make a few definitions and set some notation.
Suppose we have $H$, a compactly supported time dependent Hamiltonian
function on the symplectic manifold $(M^{2n}, \omega)$.  We may assume that 
for each $t$,
$$\min_{x \in M} H_t(x) = 0.$$ 
We write for the graph of $H$
$$\Gamma_H = \{ ( x, H_t(x), t) \} \subset M \times {\bf R} \times [0,1].$$ 
Now, let 
$$h_\infty = \max_{x \in M, t \in [0,1]} H_t(x)$$
and suppose
$\ell(t): [0,1] \rightarrow [-\delta, 0]$ is a function which is negative and
close to zero.  A thickening of the region under $\Ga_H$ is
$$R_H^-(\frac{\nu}{2}) = \{ (x,s,t) \; \vline \; \ell(t) \leq s \leq H_t(x) \}
\subset M \times [\ell(t), h_{\infty}] \times [0,1]$$
where $ \int_0^1 - \ell(t) dt = \frac{\nu}{2}$. Similarly, we can define
$R_H^+(\frac{\nu}{2})$ to be a slight thickening of the region above
$H$:
$$R_H^+(\frac{\nu}{2}) = \{ (x,s,t) \; \vline \; H_t(x) \leq s \leq \mu_H(t) \}
\subset M \times [0, \mu_H(t)] \times [0,1]$$
where $\mu_H(t)$ is a function dependent on $H$ and $t$
such that  
$$\mu_H(t) \geq \max_{x \in M} H_t(x) \mbox{ and } \int_0^1 (\mu_H(t) - h_\infty) dt = \frac{\nu}{2}.$$
We define
$$R_H(\nu) = R_H^-(\frac{\nu}{2}) \cup R_H^+(\frac{\nu}{2}) \subset
M \times {\bf R} \times [0,1].$$

For example,  consider $P$ defined on $\CP$
$$P([z_0: z_1: z_2]) = \frac{\pi}{2} \frac{|z_0|^2}{|z_0|^2 + |z_1|^2
+ |z_2|^2}. $$ 
Then,
$$\Ga(P) = \left \{ \left ([z_0: z_1: z_2] , \frac{\pi}{2} \frac{|z_1|^2}{|z_0|^2 + |z_1|^2
+|z_2|^2}, t \right ) \right \} \subset {\CP} \times [0, \frac{\pi}{2}] \times [0,1]$$
$$R_{P}^-(\frac{\nu}{2}) = \left \{\left( [z_0: z_1: z_2], s, t \right) \; \vline \; \ell(t) \leq s \leq
 \frac{\pi}{2} \frac{|z_1|^2}{|z_0|^2 + |z_1|^2+|z_2|^2} \right \}$$
$$\subset {\CP} \times [\ell(t), \frac{\pi}{2}] \times [0,1]$$
$$R_{P}^+(\frac{\nu}{2}) = \left \{\left( [z_0: z_1: z_2], s, t \right) \; \vline  \;
 \frac{\pi}{2} \frac{|z_1|^2}{|z_0|^2 + |z_1|^2+|z_2|^2}  \leq s \leq \mu_{P}(t) \right \}$$
$$\subset {\CP} \times [0, \mu_{P}(t)] \times [0,1]$$
and
$$R_{P}(\nu) = \{ \left( [z_0: z_1: z_2], s, t \right) \; \vline  \;
 \ell(t)  \leq s \leq \mu_{P}(t)\}  \subset {\CP} \times {\bf R} \times [0,1].$$

Note that we can equip $R_H^-(\frac{\nu}{2})$, $R_H^+(\frac{\nu}{2})$, and 
$R_H(\nu)$ with the product symplectic
form $\Omega = \omega \oplus ds \wedge dt$.  We need the following
definition from \cite{LM2} which describes manifolds such as 
$(R_H(\nu), \Omega)$:

\begin{definition}
Let $(M, \omega)$ be a symplectic manifold and $D$ a set diffeomorphic
to a disc in $({\bf R^2}, ds \wedge dt)$.  Then, the manifold $Q = 
(M \times D, \Omega)$ endowed with the symplectic form $\Omega$ is called
a {\bf quasi-cylinder} if 
\begin{description}
\item[(i)] $\Omega$ restricts to $\omega$ on each fibre $M \times 
\{pt\}$;
\item[(ii)] $\Omega$ is the product $\omega \oplus (ds \wedge dt)$ near the boundary
$M \times \partial D$, and, in the case where $M$ is non-compact, outside
a set of the form $X \times D$ for some compact subset $X$ in $M$.
\end{description}
\end{definition}
Note that for any Hamiltonian $H$, $(R_H(\nu), \Omega)$ is a quasi-cylinder symplectomorphic to 
$M \times D(L(H) + \nu)$ where $D(a)$ denotes the two-disk with area $a$.
Since $\Omega = \omega \oplus ds \wedge dt$ everywhere, not just
near the boundary, $R_H(\nu)$ is called a {\bf split} quasi-cylinder.  
We define the {\bf area} of any compact quasi-cylinder $(M \times D(a), \Omega)$ to
be the number  $A$ such that
$$ \mbox{ vol } (M \times D(a), \Omega) = A \cdot \mbox{ vol } (M, \omega).$$
Note that if $ (M \times D(a), \Omega)$ is split, its area is simply $a$.  
The area of $R_{P}(\nu)$, therefore, is $\nu + \frac{\pi}{2}$.  

Now, suppose $H_t$ and $K_t$ are two Hamiltonians on $M$ such that 
$\phi_1^H = \phi_1^K$ and the path $\phi_t^H$ for $ 0 \leq t \leq 1$  is
homotopic (with fixed endpoints) to the path $\phi_t^K$ in $Ham^c(M)$.  
We may join $\Gamma_K$ to $\Gamma_H$ via the map
$$g(x,s,t) = (\phi_t^H \circ (\phi_t^K)^{-1}(x), s - K(x)+H(\phi_t^H \circ {\phi_t^K}
^{-1}(x)), t).$$  This map $g$ extends to a symplectomorphism of   $R_K^+(\frac{\nu}{2})$, and we define
$$(R_{H,K}(\nu), \Omega) = R_H^-(\frac{\nu}{2}) \cup R_K^+ (\frac{\nu}{2}).$$
Because the  loop $\phi_t^H \circ {\phi_t^K}^{-1}$ is contractable in
$Ham^c(M)$, Lalonde and McDuff are able to show that $(R_{H,K}(\nu), \Omega)$
is a quasi-cylinder diffeomorphic to 
$$M \times \{s,t \in {\bf \rm R ^2} \; \vline \; \la(t) \leq s \leq \mu_H(t) \} \cong M \times D(L(H) + \nu).$$
Note that $R_{H,K}(\nu)$ is not necessarily a split quasi-cylinder, and thus the
area of $R_{H,K}(\nu)$ is not necessarily $L(H) + \nu$.

The key to the analysis in \cite{LM2} is the following lemma, whose proof we include
for the convenience of the reader.

\begin{lemma}
\label{lemma: keypoint}
\mbox { (Lalonde-McDuff,\cite{LM2}, Part II, Lemma 2.1) } Suppose that $L(K_t) < L(H_t) = A.$  Then, for sufficiently small $\nu >0$,
at least one of the quasi-cylinders $(R_{H,K}(\nu), \Omega)$ and
$(R_{K,H}(\nu), \Omega)$ has area $< A$.  
\end{lemma}

\proof{}
Choose $\nu > 0$ so that 
$$L(K_t)+2 \nu < L(H_t),$$
and suppose first that $M$ is compact.  Evidently,
$$ \begin{array}{rcl}
{\rm vol} (R_{H,K}(\nu)) + {\rm vol} (R_{K,H} (\nu)) &
= & {\rm vol} (R_H(\nu)) + {\rm vol} (R_K(\nu)) \\
& = & ({\rm vol} M) \cdot (L(H_t) + L(K_t)+2 \nu) \\
& < & 2({\rm vol} M) \cdot L(H_t) 
\end{array} $$
where $R_H(\nu) = R_H^- (\frac{\nu}{2}) \cup R_H^+ (\frac {\nu}{2})$.
If $M$ is non-compact, we may restrict to a large compact piece $X$ of $M$
and then take the volume. \QED

Lemma \ref{lemma: keypoint} tells us that if the area of both quasi-cylinders
$(R_{H,K}(\nu), \Omega)$ and $R_{K,H}(\nu), \Omega)$ is greater than or
equal to $L(H_t)$, then $L(H_t) \leq L(K_t)$.  To develop their
criteria for length minimizing paths, Lalonde and McDuff
use the theory of symplectic capacities to estimate the area of quasi-cylinders.
A symplectic capacity is a function from the set of symplectic manifolds to 
${\bf R} \cup \{ \infty \}$ satisfying certain properties; in particular, it is a symplectic 
invariant.  For more information on symplectic capacities,
see \cite{HZ}.      Suppose we have chosen a 
particular capacity $c$ and symplectic manifold $(M, \omega)$.  We say the  
{\bf capacity-area inequality} holds for $c$ on $M$ if 
$$c(M \times D(a), \Omega) \leq \mbox{ area of } (M \times D(a), \Omega)$$
holds for all quasi-cylinders $(M \times D(a), \Omega)$. In the next
section, we will give examples of manifolds and capacities that satisfy
this condition.   Although capacities are applied to symplectic manifolds,
Lalonde and McDuff define the capacity of a Hamiltonian in the following way 
\cite{LM2}.

\begin{definition} The capacity ${\bf c(H)}$ of a Hamiltonian function
$H_t$ is defined as
$$c(H) = \min \{ \inf_{\nu > 0} c(R_H^-(\frac {\nu}{2})),
\inf_{\nu > 0} c(R_H^+(\frac {\nu}{2})) \}.$$
\end{definition}  

Now, take a manifold $M$ and a capacity $c$ such that the capacity-area
inequality holds for $c$ on $M$, and suppose that we have a Hamiltonian 
$H_t : M \rightarrow {\bf R}$ for which
$$c(H) \geq L(H_t).$$
Then, for any Hamiltonian $K_t$ generating a flow $\phi_t^K$ which is homotopic
with fixed end points to $\phi_t^H$ (and thus has $\phi_1^K = \phi_1^H$), we can embed $R_H^-(\frac{\nu}{2})$ into $R_{H,K}(\nu)$ and $R_H^+(\frac{\nu}{2})$
into $R_{K,H}(\nu)$.  Thus, we know
$$ L(H_t) \leq c(H) \leq c(R_H^-(\frac{\nu}{2})) \leq c(R_{H,K}(\nu))$$
$$ L(H_t) \leq c(H) \leq c(R_H^+(\frac{\nu}{2})) \leq c(R_{K,H}(\nu)),$$
with the last inequality in both lines holding by the monotonicity property of capacities.  Since capacity-area inequality holds, we know that the
areas of both quasi-cylinders $R_{H,K}(\nu)$ and $R_{K,H}(\nu)$ must be greater than or equal to their capacities
and hence greater than or equal to $L(H_t)$.  Therefore, by Lemma 
\ref{lemma: keypoint},  $L(K_t) \geq L(H_t).$
This proves the proposition from \cite{LM2}:

\begin{prop}
\label{prop: criteria}
\mbox{ (Lalonde-McDuff, \cite{LM2}, Part II, Proposition 2.2)} Let $M$ be any 
symplectic manifold and $H_{t \in [0,1]}$ a Hamiltonian
generating an isotopy $\phi_t^H$ from the identity to $\phi_1^H$.  Suppose
there exists a capacity $c$ such that the following two conditions hold:
\begin{description}
\item[(i)] $c(H) \geq L(H_t)$ and
\item[(ii)] there exists a class ${\cal S}$ of Hamiltonian isotopies homotopic
rel endpoints to $\phi_t^H$ , $t \in [0,1]$, which is such that the 
capacity-area inequality holds (with respect to the given capacity $c$) for all
quasi-cylinders $R_{H,K}(\nu)$ and $R_{K,H}(\nu)$ corresponding to Hamiltonians
$K_t \in {\cal S}.$  
\end{description}
Then, the length of the path $\phi_t^H$ is minimal among all paths in 
${\cal S}.$
\end{prop}

Hence, to show that $H_t$ generates a length minimizing
 path $\phi_t^H$ for $t \in [0,1]$ among all paths homotopic rel endpoints,
we need only produce a capacity $c$ that satisfies the above conditions (i)
and (ii).  In
fact, Lalonde and McDuff show that if the capacity-area inequality holds
for all split quasi-cylinders of the form $M \times D(a)$, then
it also holds for all $R_{H,K}$ in Proposition 4.4 of \cite{LM2}.  Therefore, it will be enough
to find a capacity that satisfies (i) and  satisfies (ii) for all split quasi-cylinders,
$M \times D(a)$.  Our ${\cal S}$ will be the set of all Hamiltonians $K_t$ where 
$\phi_1^K = \phi_1^H$ and $\phi_t^K$ is homotopic rel endpoints to $\phi_t^H$.

\subsection{Capacities}

The symplectic capacities we will work with in this paper are the Gromov capacity,
$c_G$,  and the Hofer-Zehnder capacity, $c_{HZ}$. We recall their definitions
for the convenience of the reader.

\begin{definition} Let $(N, \omega)$ be a symplectic manifold of dimension
$2n$. 
\begin{description}
\item[(i)] The {\bf Gromov capacity} 
$$c_G(N, \omega)=
\sup \left \{ \pi r^2 \; \vline \begin{array}{c}
  \exists \mbox{ a symplectic embedding } \\
\phi: (B^{2n}(r), \omega_0) \rightarrow (N^{2n}, \omega) \end{array}
\right \}$$ 
where $(B^{2n}(r), \omega_0)$ is the open $2n$-dimensional ball
with radius $r$ endowed with the standard symplectic form. 
\item[(ii)] The {\bf Hofer-Zehnder capacity} 

$${\cHZ}(N, \omega) = \sup \{ \max(H) \; \vline \; H \in {\cal H}_{ad}(N, \omega) \}$$
where $ {\cal H}_{ad}(N, \omega)$ consists of all of the autonomous Hamiltonians on $N$ 
satisfying the properties
\begin{description}
\item[(a)] There exists a compact set $\kappa \subset N \setminus \partial N$
depending on $H$ so that \linebreak
$H \; \vline \; (N \setminus \kappa) = \max(H)$ is constant.
\item[(b)] There is a nonempty open set $U$ depending on $H$ such that 
$H \; \vline \; U = 0$.
\item[(c)] $0 \leq H(x) \leq \max(H)$ for all $x \in N$.
\item[(d)] All $T$-periodic solutions of the Hamiltonian system $\dot x =
X_H(x)$ on $N$ with $0 \leq T \leq 1$ are constant.
\end{description}
\end{description}

\end{definition}  

To check that the capacity-area 
inequality holds on split quasi-cylinders for either of these capacities is a non-trivial procedure.  By
using J-holomorphic curve techniques, Lalonde and McDuff show in \cite{LM2}
that it
holds for $c_G$ on manifolds $M$, compact at $\infty$, which are of 4 dimensions or
fewer or which are semi-monotone.  Recently, they have shown 
that it holds for all $M$ in \cite{LM1} .

We know, then, that  condition (ii) from
Proposition \ref{prop: criteria} is satisfied for $c_G$ on any manifold,
and in particular on $\CP$ endowed with the standard symplectic form
$\tau_0$ derived from the Fubini-Study metric .
In the appendix, we use Propostion \ref{prop: criteria} and $c_G$,
construct specific embeddings of 6-balls, and  show that rotation 
through $\pi$ radians around the first coordinate 
in $\CP$ and in $\CPb$ (the blow-up of $\CP$ at the point $[1: 0: 0]$)
is  length minimizing among all homotopic paths.  In addition we explain why 
$c_G$, for volume reasons, cannot be used to show the analagous rotation around the second 
coordinate in $\CPb$ is length minimizing.    

Since $c_{HZ}$ is not directly related to volume in the same way as
$c_G$, the next natural step is to see if we can use $c_{HZ}$ to show paths,
in particular this rotation in the second coordinate of $\CPb$, are
length minimizing.  Thus
we need to examine the conditions under which the capacity-area inequality (condition (ii) of 
Proposition \ref{prop: criteria})
holds for $\cHZ$.  Recall that a symplectic manifold $(M, \omega)$ is 
{\bf weakly exact} if $\omega$ restricted to $\pi_2(M)$ is zero.  The following theorem 
from \cite{HV} is quoted as Theorem 1.17 in \cite{LM2}:

\begin{theorem}
\label{theorem: caHV}
(Hofer-Viterbo)  Let $(M, \omega)$ be a compact symplectic manifold which
is weakly exact. 
Then for all $a>0$,
$$\cHZ (M \times D(a), \omega \oplus \sigma) \leq a.$$  
\end{theorem}

However, $\CPb$ is not weakly exact, as the Hurewicz homomorphism is
an isomorphism  between $H_2(\CPb, {\bf Z})$ and $\pi_2(\CPb)$.  In order
to eventually apply Proposition \ref{prop: criteria} to $\CPb$ using $\cHZ$,
we will go back to the original proof of  Theorem \ref{theorem: caHV}
and show that the restriction that  $M$ is weakly exact can be changed to 
$M$ has dimension two or four.   Hence, in the next section we arrive at 

\begin{theorem}
\label{theorem: chzprod}
 Let $(M, \omega)$ be a compact symplectic manifold of dimension two
or four.   Then for all $a >0$, 
$$\cHZ (M \times D(a), \omega \oplus \sigma) \leq a.$$
\end{theorem}

Theorem \ref{theorem: chzprod} enables us to prove the following main result.
\begin{theorem}
\label{theorem: mainthmgeod}
Let $(M, \omega)$ be a compact symplectic manifold of dimension two or four. 
Let $\phi^H_t$ for
$0 \leq t \leq 1$ be a path in $Ham (M)$  generated by an autonomous 
Hamiltonian 
$H: M \rightarrow {\bf R}$  such that $\phi^H_0$ is the identity
diffeomorphism and  $\phi^H_t$ has no non-constant closed trajectory in
time less than 1.
Then, $\phi^H_t$ for $0 \leq t \leq 1$ is length minimizing among all 
homotopic paths between the identity and $\phi^H_1$.
\end{theorem}

Finally,  a consequence of Theorem \ref{theorem: mainthmgeod}
and Proposition \ref{proposition: r1}
 is that the path $\phi_t$ for $0 \leq t \leq 1$ in $Ham(\CPb)$ given by
$$\phi_t [z_0: z_1: z_2] = [z_0: e^{\pi it}z_1: z_2]$$
is length minimizing between the identity ($\phi_0$) and rotation
by $\pi$ radians in the second coordinate ($\phi_1$).    
 
\section{The capacity-area inequality for $c_{HZ}$}
In the first part of this section, we analyze the proof of Theorem
\ref{theorem: caHV} which states sufficient conditions on $M$ for
$c_{HZ}$ to satisfy the capacity-area inequality on $M$.  Then, in
the second portion, we show that the weakly exact hypothesis in this
theorem can be changed to dimension two or four.

\subsection{Hofer and Viterbo's proof of Theorem 2.6}
We now examine Hofer and Viterbo's proof of Theorem \ref{theorem: caHV}
to determine why they need the weakly exact condition \cite{HV}.  Unfortunately, 
their notation is different from the notation in \cite{LM2}, so we will first need to 
provide some sort of dictionary to explain the theorem as they have stated
it.  

Let $[S^2, M]$ be the set of homotopy classes of maps from $S^2$ to $M$.
We apply $\omega$ to such a class $\alpha \in [S^2, M]$ by evaluating $\omega$ 
on the representative of $\alpha$ in $H_2(M, {\bf Z})$. Define
$$m(M, \omega) = \inf \{ \langle \omega, \alpha \rangle \; \vline \; \alpha \in
[S^2, V], \; 0 < \; \langle \omega, \alpha \rangle \}.$$

Note that if $M$ is weakly exact, $m(M, \omega) = \infty$.    If
for some particular class $\alpha \in H_2(M)$ we have $\langle \omega, \alpha
\rangle = m(M, \omega)$, then $\alpha$ is called $\omega$-minimal.   The
theorem of Hofer and Viterbo which is equivalent to Theorem \ref{theorem: caHV} is

\begin{theorem}
\label{theorem: HV1.12}
(Hofer-Viterbo, \cite{HV}, Theorem 1.12) Let $(M, \omega)$ be a compact symplectic manifold and let $\sigma$ be
a volume form for $S^2$ such that 
$\int_{S^2}\sigma = a$  and 
$$0 < a \leq m(M, \omega).$$
Suppose $K: M \times S^2(a) \rightarrow {\bf R}$ is a smooth (time
independent) Hamiltonian
such that
$$K \; |_{{\cal U}(\ast)} = k_0 \mbox{ and } K \; |_{{\cal U}(M \times \{ \infty \})} = k_{\infty}$$
for suitable neighborhoods of $M \times \{ \infty \}$ and some point 
$ \ast \not\in M \times \{ \infty \}$.  Suppose 
$$\begin{array}{ccc}
k_0 < k_{\infty} & and & k_0 \leq K \leq k_{\infty}.
\end{array}$$
Then, the Hamiltonian system $\dot x = X_K(x)$ on the symplectic manifold
\linebreak
$(M \times S^2(a), \omega \oplus \sigma)$ possesses a non-constant
$T$-periodic solution with
$$0< (k_{\infty} - k_0)T  < a.$$
\end{theorem}

The task now at hand is to see why Theorem \ref{theorem: HV1.12} is equivalent to Theorem \ref{theorem: caHV}.  Remember that 
$${\cHZ}(N, \omega) = \sup \{ \max(H) \; \vline \; H \in {\cal H}_{ad}(N, \omega) \}
$$
where $ {\cal H}_{ad}(N, \omega)$ consists of all of the autonomous Hamiltonians on $N$ 
satisfying the properties:
\begin{description}
\item[(a)] There exists a compact set $\kappa \subset N \setminus \partial N$
depending on $H$ so that $H \; \vline \; (N \setminus \kappa) = \max(H)$ is constant.
\item[(b)] There is a nonempty open set $U$ depending on $H$ such that 
$H \; \vline \; U = 0$.
\item[(c)] $0 \leq H(x) \leq \max(H)$ for all $x \in N$.
\item[(d)] All $T$-periodic solutions of the Hamiltonian system $\dot x =
X_H(x)$ on $N$ with $0 \leq T \leq 1$ are constant.
\end{description}

Clearly, proving Theorem \ref{theorem: caHV}  is the same as showing that
any properly normalized Hamiltonian $K$ on $M \times D(a)$ with 
$\max(K) > a$ 
has a non-constant orbit with period $T \leq 1$.  In Theorem \ref{theorem: HV1.12}, 
Hofer and Viterbo consider the completion $M \times S^2(a)$ of $M \times D(a)$.  For 
simplicity, we will also denote the symplectic form on
$S^2(a)$ by $\sigma$.  The neighborhood $U(M \times \infty) 
\subset M \times S^2(a)$ corresponds to a neighborhood of $\partial(M \times D(a))$ in
Theorem \ref{theorem: caHV}.  The hypotheses concerning
the values $k_0$ and $k_\infty$ in Theorem \ref{theorem: HV1.12} correspond to the 
conditions (a) 
(b), and (c) describing the requirements for $K$ to be a member of ${\cal H}_{ad}$.
The hypothesis $0 < a \leq m(M, \omega)$ in Theorem 
\ref{theorem: HV1.12} is satisfied for all $a$ if and only if $M$ is weakly
exact.   Finally, the quantity $k_{\infty} - k_0$ corresponds to $\max(K)$.  Hence, to
show the equivalence of the two theorems
we need to suppose in  Theorem \ref{theorem: HV1.12} that $k_{\infty} - k_0 \geq a$ and
show that we get a closed non-constant orbit of period  $T \leq 1$. In fact, the conclusion of Theorem \ref{theorem: HV1.12}
tells us exactly that  we get a non-constant orbit of period
$$ T <  \frac{a}{k_{\infty} - k_0},$$
so that if $k_{\infty} - k_0 \geq a$ then $T \leq 1$.

We eventually want to prove Theorem \ref{theorem: HV1.12} without the hypothesis
$a \leq m(M, \omega)$.  Consider the symplectic manifold $(M \times S^2(a), \omega 
\oplus \sigma)$ and let
${\cal J}$ be the set of all smooth almost complex structures $J$ compatible 
with $\omega \oplus \sigma$ 
on $M \times S^2(a)$.  The original proof of Theorem \ref{theorem: HV1.12}
uses $J-$holomorphic curves with a split compatible
almost complex structure $J \in {\cal J}$ on $M \times S^2(a)$  that is
regular for the class $A = [\{pt\} \times S^2]$ in the sense of Theorem 3.1.2 of
\cite{MS}.  Hofer and Viterbo use a split $J$ so that they
can easily verify the condition $a < m(M, \omega)$ in certain settings.
Since this condition is exactly the hypothesis we will remove, in this
discussion we do not need to restrict ourselves to a split $J$.  We will,
however, need to impose more regularity conditions on $J$ later. 

After a $J$ is fixed,
the proof of Theorem \ref{theorem: HV1.12} proceeds by determining the
$S^1$-cobordism class of a certain moduli space of $J$-holomorphic spheres
whose image is in $M \times S^2(a)$.  This moduli space  ${\cal M}(J)$
consists of the set of maps $u \in C^\infty(S^2, M \times S^2(a))$ that
satisfy

$$[u] = [\{ {\rm pt} \} \times S^2(a)] = A \in H_2(M \times S^2(a), {\rm \bf Z})$$
$$ \int_D u^\ast (\omega \oplus \sigma) = \frac{a}{2} \; {\rm  where} \;
D = \{ z \; \vline \; |z| \leq 1 \}$$
$$u(0) = \{ \ast \}, \; \; u(\infty) \in M \times \{ \infty \}$$
$$\overline{\partial}_Ju = 0.$$

Hofer and Viterbo show the $S^1$-cobordism class of ${\cal M}(J)$ is not zero and hence a related family
${\cal C}$ of perturbed $J$-holomorphic spheres is not compact.  
Specifically, 
$${\cal C} = \{ (\lambda, u) \in [0, \infty) \times {\cal B} \; \vline \;
 \bar \partial_J u + \lambda k(u) = 0 \}$$
where $k(u)$ is basically a scaling of the gradient of $K$  and 
${\cal B}$ is the set of maps $ u \in H^{2,2}(S^2, M \times S^2(a))$ 
that satisfy 
$$[u] = A \in H_2(M \times S^2(a), {\rm \bf Z})$$
$$ \int_D u^\ast (\omega \oplus \sigma) = \frac{a}{2} \; {\rm  where} \;
D = \{ z \; \vline \; |z| \leq 1 \}$$
$$u(0) = \{ \ast \}, \; \; u(\infty) \in M \times \{ \infty \}.$$

We can see that given a $\lambda$, the map $u$ for $(\lambda, u) \in {\Cc}$
is almost fixed. Since $J$ is regular, the dimension of the moduli space of 
perturbed $J-$holomorphic spheres of class $A$
is $2c_1(A) + {\rm dim}(M) + 2 = 6 + {\rm dim}(M)$ (\cite
{MS}, Theorem 3.12).  However, ${\cal C}$ does not consist of all
of these spheres; the restrictions placed upon the elements
in ${\cal B}$ reduce the dimension of ${\cal C}$ greatly. The first
normalization condition on the area imposes a loss of 1 dimension.
The next restriction, fixing the image of $\{0\}$, imposes a loss
of ${\rm dim}(M) + 2$ dimensions.  Finally, restriction the image
of $\{\infty\}$ results in a loss of 2 dimensions.  Hence, the set
of spheres we are considering in the second factor of ${\cal C}$
will have dimension $6 + {\rm dim}(M) -
1-({\rm dim}(M) + 2) - 2 = 1$.  This degree of freedom corresponds to rotation by $S^1$ 
of $S^2$.  We are basically fixing the parametrization of $u$ except for allowing this $S^1$ action.
Note, then, that ${\Cc}$ is a two dimensional space: one dimension for
the $\lambda$ coordinate and one dimension which corresponds to this
$S^1$ rotation.

Hofer and Viterbo analyze the noncompactness of ${\Cc}$ and show that it cannot be due to
a bubbling off of perturbed $J$-holomorphic curves. Since
there are no bubbles, there are uniform bounds on the derivatives of the 
$u$.  They view the $u$ not as maps from the sphere, but rather as maps
from the non-compact cylinder $S^1 \times {\bf R}$.  Hence, $\cal C$ consists
 of maps 
with finite energy whose domain is 
an infinitely long cylinder.  In the same manner as in Floer theory,
Hofer and Viterbo show the noncompactness of
${\cal C}$  produces a sequence of maps that converge to a 
closed non-constant orbit $x$  which is a solution of the equation 
$\dot x = X_K(x)$.  

When we remove the restriction 
$a \leq m(M, \omega)$, each of the steps in the
proof of Theorem \ref{theorem: HV1.12} goes through with only minor 
adjustments, except for the proof
of the statement that there are no bubbles.  It turns out, however, that
this difficulty can be overcome.  In the next section, we
give a new proof that shows that it is
still true generically that
no sequence of elements in ${\cal C}$ converges to a bubble when we remove the area
restriction if $M$ has dimension two or four.  We need the dimension restriction
on $M$ to rule out the possibility of multiply covered curves on $M \times S^2$.  

\subsection{Noncompactness in ${\cal C}$ cannot be due to bubbling}
We will show that for generic $J \in {\cal J}$, the space of bubbles 
which are limits of sequences of elements in ${\cal C}$ is empty.  We first show
that for generic $J$, the space of cusp curves which have two components
in empty.  

There are five distinct types of two-component bubbles which are possible.  
We consider them separately.    For each type, we will find a dense set of $J$ so that 
the particular type does not occur; the intersection of these five sets will be our set 
of regular $J$.  The first case is when the point $z_0$ where the derivitave blows up in 
$S^2$ lies on the upper hemisphere but is not $\{\infty \}$; the other
cases are when the point lies on the lower hemisphere, when the point lies on
the equator, when the point is $\{0\}$, and when the point is  $\{\infty\}$.   We must
separate the cases in this way to handle  appropriately the restrictions of  curves
that lie in $\Cc$.

Let us now investigate the first case.  We will represent the $\lambda k(u)$ perturbed 
component of the cusp curve
by the class $A-Y$ and the $J$-holomorphic bubble by the class $X$.  Let
us for now assume that $X=Y$, and therefore that the homological sum of the
two component classes is $A$. Note that this need not be the case: since
we only consider simple cusp curves as limiting elements, we may have had to 
reduce a multiply covered curve and thus have lost some homology.  We will
discuss this later on and see that, since  $M$ has dimension two or 
four, it poses no obstacle.

 Define the universal moduli spaces 
$$\mu^{\lambda}(A-Y, {\cal J}) = \{(u, J) \; \vline \; u: S^2 \rightarrow M \times S^2(a),
[{\rm Im}(u)] = A-Y, \; \bar \partial_J u + \lambda k(u) = 0\}$$
and 
$$\mu(Y, {\cal J}) = \{(v, J) \; \vline \; v: S^2 \rightarrow M \times S^2(a),
[{\rm Im}(v)] = Y, \; \bar \partial_J v = 0 \}.$$
We will write $\mu^{\lambda}(A-Y, J)$ or $\mu(Y,  J)$ when we wish to
refer to the moduli space consisting of curves corresponding to a single $J$.

\begin{figure}[htb]
\centering
\includegraphics{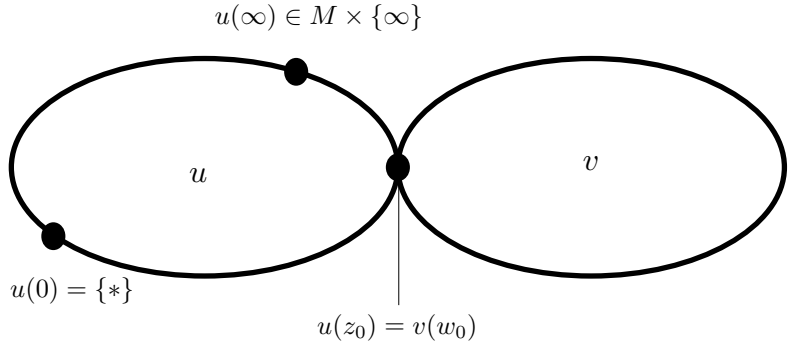}
\caption{Bubbling in the upper hemisphere}
\label{figbubble}
\end{figure}

We must show that for a generic $J$, the subset  of elements in \linebreak
$\mu^{\lambda}(A-Y, J) \times  \mu(Y,  J)$ which could be a limit of curves satisfying 
the restrictions
 of ${\cal C}$ and which are bubbles is empty.  We are considering the first type
of bubble where the point at which the derivative blows up to form the bubble lies
in the upper hemisphere.  A picture 
of the cusp curve is shown in Figure \ref{figbubble}.

Let ${\cal U}^\lambda$ be the space 
$${\cal U}^\lambda =  \bigcup_{J \in \Jj} 
\mu^{\lambda}(A-Y,J) \times \mbox{ upper hemisphere of }S^2 \times \mu(Y,J) \times_G S^2. $$
Here $G= {\rm PSL}(2, {\bf C})$ is the six dimensional holomorphic 
reparametrization group of $S^2.$   (Note that for different types of bubbles we will be able to
quotient by different symmetry groups.)  Define the space 
$$\Uu = \{\lambda, \Uu^\lambda \; \vline \; \lambda \in [0, \lambda_\infty]\}$$
where $\lambda_\infty$ is a constant described in \cite{HV} that depends on $A$ and $M$.  
We let $\Uu_J$ be the restriction of $\Uu$ to a particular $J \in \Jj$.  Next, we isolate the
curves in $\Uu$ which are bubbles and can be limits of a sequence of elements in
$\Cc$.

Consider the evaluation map $ev$ where
$$ev: \Uu \rightarrow (M \times S^2) ^4 \times {\bf R}$$
by
$$ ev(\lambda, J,u,z_0, v,w_0) = \left( u(\infty), u(0), u(z_0), v(w_0), 
\int_D u^*(\omega \oplus \sigma) \right).$$
Let
$${\cal D} = ev ^{-1} \left( (M \times \{\infty\} ), \{*\}, \Delta,   \frac{a}{2} \right)$$
where $\Delta$ stands for the diagonal in $(M \times S^2) \times 
(M \times S^2)$. We let ${\cal D}_J$ be the restriction of ${\cal D}$ to a particular $J \in \Jj$ .  
 Note that $\Dd_J$  consists exactly of the elements that are
bubbles and could be the limit of a sequence of elements in $\Cc$.
Our aim is to prove that, for a generic $J$, $\Dd_J$ is empty.   
Recall that although the curves in $\Cc$ do not have a full reparametization group acting
on them, they do have an $S^1$ action.  
Hence, if $\Dd_J$ is non-empty, it must be of dimension at least one to account for this
symmetry.  We will show that generically $\Dd_J$ has dimension zero, and therefore 
it must be empty.  

Of course for the analysis to make sense in this infinite dimensional setting, $\Uu$ needs
to be a Banach manifold.  Hence we must restrict the set of almost complex structures
$\Jj$ to contain only those with sufficient smoothness and require that the curves $u$
and $v$ belong to an appropriate Sobolev space.  These specific notions are described
explicitly  in Propositions 6.2.2 and 3.4.1 of \cite{MS}.  

\begin{proposition}
\label{proposition: DJempty}
There exists a set of complex structures $\Jj_0 \subset \Jj$ of second 
category such that for $J \in \Jj_0$, $\Dd_J$ is empty.
\end{proposition}

\proof{}
We follow the steps in the proof of Theorem 6.3.2 from \cite{MS}, using three
lemmas.  
 The first shows that $ev$ is transversal  onto the set 
$$  \left( (M \times \{\infty\} ), \{*\}, \Delta,  \frac{a}{2} \right).$$ 
Then, 
the second, proves that the
projection map $\pi$ from $\Uu$ onto $\Jj$  is a Fredholm operator.  Hence, there will be a set
$\Jj_{reg}$ of regular values of $\pi$  such 
that for $J\in\Jj_{reg}$,  $\Uu_J$ is a manifold.  We prove that the point 
$(\lambda, J, u,z_0,v,w_0)$  is a regular point of $\pi$ exactly when $ev$ restricted 
to $\Uu_J$ is transversal to $\left( (M \times \{\infty\} ), \{*\}, \Delta,   \frac{a}{2} \right)$
at the point $ev(\lambda, J, u,z_0,v,w_0)$ .  
Therefore, for  $J \in \Jj_{reg}$, the submanifold $\Dd_J \subset \Uu_J$  has
its expected codimension.  The third lemma states that for $J$ in some subset $\Jj_0 \subset \Jj_{reg}$ 
this codimension is equal to the dimension of
$\Uu_J$ and hence $\Dd_J$ is empty.  

\begin{lemma}
\label{lemma: evtrans}
The map $ev$ is transversal to 
$$ \left( (M \times \{\infty\} ), \{*\}, \Delta, \frac{a}{2} \right)\subseteq (M \times S^2)^4 \times 
{\bf R}.$$ 
\end{lemma}

\proof{}
Define the map 
$$ev^{z_0, w_0} :  \bigcup_{\lambda \in [ 0, \lambda_{\infty} ]}  \Jj \times \mu^\lambda(A-Y,J)
\times \mu(Y,J) \rightarrow (M \times S^2)^4 \times {\bf R}$$
by
$$ev^{z_0,w_0} ( \lambda, J,u,v) = (u( \{ \infty \}), u(0), u(z_0), v(w_0), 
\int_D u^*(\omega \oplus \sigma)).$$
It suffices to show that for some pair 
$$(z_0, w_0) \in ( \mbox{upper hemisphere of }S^2  - \{ \infty \}  \times S^2),$$
the map $ev^{z_0, w_0}$ is transversal onto  
$$ \left( (M \times \{\infty\} ), \{*\}, \Delta,  \frac{a}{2} \right)\subseteq (M \times S^2)^4 \times 
{\bf R}.$$ 
 
Let $\pi_i : (M \times S^2)^4 \times {\bf R} \rightarrow (M \times S^2)$ be projection onto the
$i$th $M \times S^2$ factor and let $\rho: (M \times S^2)^4 \times {\bf R} \rightarrow {\bf R}$
denote projection onto the last factor.  
Since transversality is a local condition and the points 
$\{ 0 \}, \{ \infty \}, $ and $z_0$
are separated, it is enough to show that 
$$ e_1 = \pi_1 \circ ev^{z_0, w_0} \mbox{ is transversal to }  M \times \{\infty\} \subset
M \times S^2 $$
$$ e_2 = \pi_2 \circ ev^{z_0, w_0} \mbox{ is transversal to } \{*\} \in M \times S^2$$
$$ e_{34} = \pi_3 \times \pi_4 \circ ev^{z_0, w_0} \mbox{ is transversal to } \Delta \subset
(M \times S^2) ^2 $$
and
$$e_5 = \rho \circ ev^{z_0, w_0} \mbox{ is transversal to } \frac{a}{2}.$$

We now recall a theorem from \cite{MS} (Theorem 6.1.1). Let $B$ be any class in
$H_2(M \times S^2).$  For $x_0 \in S^2$, they define the map
$$e^{x_0}: \mu (B, \Jj) \rightarrow M \times S^2 \mbox{ by } e^{x_0}(u,J) = u(x_0).$$

\begin{theorem}
\label{theorem: submersion}
(McDuff-Salamon)  For any point $x_0$, the map $e^{x_0}$ is a
submersion onto $M \times S^2$.
\end{theorem}
This theorem is stated for unperturbed curves, but its proof applies to the perturbed
case as well.  Hence,
 it directly implies that $e_1$ and $e_2$  are submersions and therefore certainly transversal.  
To show that $e_{34}$ is 
transversal to $\Delta$, note that the normal
bundle to $\Delta$ at the point $(q,q) \in (M \times S^2)^2$ is spanned by 
$0 \oplus T_q(M \times S^2)$.  By applying Theorem \ref{theorem: submersion}
 we see $e_4 =\pi_4 \circ ev^{z_0, w_0}$ is a submersion.  Therefore, $e_{34}$ is indeed transversal to
$\Delta$.  Finally, $e_5$ is transversal to $\frac{a}{2}$ simply because the area over $D$ 
of the pull back of the symplectic form by the curve $u$  has neither 
a local maximum nor a local minimum at $\frac{a}{2}$. \QED

\begin{lemma}
\label{lemma: codimDJ}
There exists a set of second category $\Jj_{reg} \subset \Jj$, so
that the codimension of $\Dd_J$ in $\Uu_J$ is $4n+7$ for all $J \in \Jj_{reg}$.
\end{lemma}

\proof{}
Consider the projection map $\pi : \Uu \rightarrow \Jj$.  Note that $(\pi)_*$ is onto,
so its cokernel is $0$ and hence finite dimensional.
At the point $(\lambda, J, u, z_0, v, w_0) \in \Uu$, the kernel of $(\pi)_*$ consists
of vectors $$\{ (\hat{\lambda}, Z, \xi_u, \hat{z}, \xi_v, \hat{w}) \; \vline \; Z = 0\}.$$
The $\xi_u$ and $\xi_v$ directions contribute finitely many dimensions no matter
which $J$ is chosen, and the other directions contribute a total of five dimensions.
Hence the kernel is finite dimensional and $\pi$ is Fredholm.

Let $\Jj _{reg} \subset \Jj$ denote the set of almost complex structures which are
regular values of $\pi$.  It is a set of second category.  
For $J \in \Jj _{reg}$, we know $\Uu _J$ is a manifold of the expected dimension. 

Now, we show that the point $(\lambda, J, u, z_0, v, w_0)$ is a regular point of $\pi$ 
if and only if at this point the restricted evaluation map
$$ev: \Uu _J \rightarrow (M \times S^2)^4 \times {\bf R}$$ 
is transversal to the set   
$\left( (M \times \{\infty\} ), \{*\}, \Delta,   \frac{a}{2} \right).$

We already know from Lemma \ref{lemma: evtrans} that the set of vectors
$$ev_*(\lambda, J, u, z_0, v, w_0)(\hat{\lambda}, Z, \xi_u, \hat{z}, \xi_v, \hat{w})$$
is transversal to $\left( (M \times \{\infty\} ), \{*\}, \Delta,  \frac{a}{2} \right),$ and
now we must explain why the subset of these vectors with $Z =0$ (corresponding to
keeping $J$ constant) is still transversal.   By Lemma \ref{lemma: evtrans} and the 
linearity of $ev_*$ we see
$$ \begin{array} {rcl}
T_{ev(\lambda, J, u, z_0, v, w_0)}((M \times S^2)^4 \times {\bf R}) & = &  \mbox {Span}
\left( \mbox{Im}(ev_*) +
T_{ev(\lambda, J, u, z_0, v, w_0)} \left( (M \times \{\infty\} ), \{*\}, \Delta,  \frac{a}{2} \right) 
\right)\\
 & = & \mbox{Span} ( \mbox{Im}(ev_* \vline_{Z=0}) + \mbox{Im}(ev_* \vline_S) + \\
& & 
T_{ev(\lambda, J, u, z_0, v, w_0)} \left( (M \times \{\infty\} ), \{*\}, \Delta,  \frac{a}{2} \right)).
\end{array}$$
Here $S$ is the set of tangent vectors that satisfy
 $\hat{\lambda} = \xi_u = \hat{z} = \xi_v = \hat{w} = 0$; that is the set of vectors for
which all of the components except possibly the one in the $Z$ direction are zero.   
However, we claim that $(\lambda, J, u, z_0, v, w_0)$ is a regular point of $\pi$ if and
only if
$$\mbox{Span} \left( \mbox{Im}(ev_* \vline_S) +
T_{ev(\lambda, J, u, z_0, v, w_0)} \left( (M \times \{\infty\} ), \{*\}, \Delta, \frac{a}{2} \right)
\right) =$$
$$\mbox{Span}\left( T_{ev(\lambda, J, u, z_0, v, w_0)} \left( (M \times \{\infty\} ), \{*\}, 
\Delta,  \frac{a}{2} \right) \right).$$
Therefore, at a regular point of $\pi$,    
$$ \begin{array} {rcl}
T_{ev(\lambda, J, u, z_0, v, w_0)}((M \times S^2)^4 \times {\bf R}) =   
\mbox{Span} \left( \mbox{Im}(ev_* \vline_{Z=0}) +   
T_{ev(\lambda, J, u, z_0, v, w_0)} \left( (M \times \{\infty\} ), \{*\}, \Delta,  \frac{a}{2} \right)
\right) \end{array}$$
and $ev$ restricted to $\Uu_J$  is transversal.
To prove the claim, note that for $(\lambda, J, u, z_0, v, w_0)$ to be a regular point of $\pi$
means that  for any $Z_0 \in T_J \Jj$, there exists a tangent vector 
$$(\hat{\lambda}, Z_0, \xi_u, \hat{z}, \xi_v \hat{w}) \in T_{(\lambda, J, u, z_0, v, w_0)} \Uu .$$
In other words, for any $Z_0 \in T_J \Jj$,
$$ev_*(\lambda, J, u, z_0, v, w_0)(\hat{\lambda}, Z_0, \xi_u, \hat{z}, \xi_v \hat{w}) \in
T_{ev(\lambda, J, u, z_0, v, w_0)}\left( (M \times \{\infty\} ), \{*\}, \Delta,  \frac{a}{2} \right).$$
Thus, adding the vectors $\mbox{Im}(ev_* \vline_S)$ to the set 
$$T_{ev(\lambda, J, u, z_0, v, w_0)}\left( (M \times \{\infty\} ), \{*\}, \Delta,  \frac{a}{2} \right)$$
do not change its span, because
these vectors only have $Z$ components and all $Z$ components are already accounted for
in the set.

Hence, for $J \in \Jj_{reg}$,  $ev$ restricted to $\Uu _J$ is transversal to 
$ ( (M \times \{\infty\} ), \{*\}, \Delta,  \frac{a}{2} )$.  The 
inverse image of this set under $ev$, called $\Dd_J$, will be a manifold of the same codimension.

$$\begin{array}{rcl}
\mbox{ codimension of } \Dd_J & = & \mbox{ codimension of } 
                                    ( (M \times \{\infty\} ), \{*\}, \Delta,  \frac{a}{2}  ) \\
                              & = & 2+(2n+2)+(2n+2) +1\\
                              & = & 4n + 7.
\end{array} $$ \QED

Finally, we calculate the dimension of $\Uu_J$. 
\begin{lemma}
\label{lemma: dimUJ}
There exists a set of second category $\Jj'_{reg} \subset \Jj$
so that for $J \in \Jj'_{reg}$, the dimension of $\Uu_J$ is $4n+7$.
\end{lemma}

\proof{}
We recall a theorem from \cite{MS} [Theorem 3.1.2].  Let $B$ be any
class in \linebreak $H_2(M \times S^2)$.
\begin{theorem}
(McDuff - Salamon) There exists a set of second category $\Jj'_{reg}(B) \subset \Jj$, 
such that for $J \in \Jj'_{reg}(B)$ the moduli space 
$\mu(B, J)$ is a smooth manifold of dimension $2c_1(B) + 2n + 2$.
\end{theorem} 
If we let our classes be $A-Y$ and $Y$, we see that
for
$$J \in \Jj'_{reg}(A-Y) \cap \Jj'_{reg}(Y) = \Jj'_{reg},$$
we have
$$\begin{array} {rcl}
\dim {\Uu _J} & = & 1 + (2c_1(A-Y) + 2n+2) + 2 + (2c_1(Y)+ 2n + 2) + 2 - 6 \\
           & = &  2c_1(A) + 4n +3 \\
           & = & 4n+7.
\end{array} $$ \QED

Let $\Jj_0 = \Jj_{reg} \cap \Jj'_{reg}$.
For  $J \in  \Jj_0$, $\Uu_J$ is a manifold of dimension 
$4n + 7$  in which $\Dd_J$ has codimension
$4n + 7$. Hence, for these $J$, $\Dd_J$ will  have dimension zero and in fact
be empty as described earlier.  Note that $\Jj_0$ is  of second category since 
it is the intersection of  two second category sets . Thus, we have proven 
Proposition \ref{proposition: DJempty}. 
\QED

\begin{prop}
\label{prop: nobub}
Suppose  $(M, \omega)$ is a compact symplectic manifold of dimension two or four. Then, there
exists a  set of second category of regular almost complex structures on 
$M \times S^2$ for which the space of bubbles which are limits of sequences of 
elements in ${\Cc}$ will be empty. 
\end{prop}

\proof{}
Proposition \ref{proposition: DJempty} tells us that for generic $J$, the space of such
bubbles that are cusp curves with two components, neither of which is multiply 
covered, where the bubble is of a given homology class $Y$ and is formed by the 
derivative blowing up at a point on the upper hemisphere,
is empty.  To deal with other types of bubbling in a two component cusp
curve is similar.  We must be careful in defining the evaluation map to 
use the correct domain, quotienting out by the appropriate reparametrization group,
and set the area condition of the last component properly.  Here are the precise 
variations, indexed by the point  on the sphere at which the derivative blows up:

\begin{description}
\item[(i) Lower Hemisphere]    Change the domain of $ev$ by setting
$${\cal U}^\lambda =  \bigcup_{J \in \Jj} 
\mu^{\lambda}(A-Y,J) \times \mbox{ lower hemisphere of }S^2 \times \mu(Y,J) \times_G S^2. $$
Use
$$\int_D u^*(\omega \oplus \sigma) + \int_{S^2} v^*(\omega \oplus \sigma)$$
for the final component in the map $ev$.
\item[(ii) Equator] Change the domain of $ev$ by setting
 $${\cal U}^\lambda =  \bigcup_{J \in \Jj} 
\mu^{\lambda}(A-Y,J) \times \mbox{ equator of }S^2 \times \mu(Y,J) \times_G S^2. $$
Use
$$\int_D u^*(\omega \oplus \sigma) + \int_{D} v^*(\omega \oplus \sigma)$$
for the final component in the map $ev$.
\item[(iii) The point $\{0\}$] Change the domain of $ev$ by setting
$${\cal U}^\lambda =  \bigcup_{J \in \Jj} 
\mu^{\lambda}(A-Y,J) \times \mu(Y,J) \times_{G_0} S^2 $$
where $G_0$ is the four dimensional set of holomorphic maps from $S^2$ to itself which
fix $\{0\}$, and let
$$ev(\lambda, J, u,v,w_0) = \left( u(\infty), v(0), u(0), v(w_0),  
\int_D u^*(\omega \oplus \sigma) + \int_{S^2} v^*(\omega \oplus \sigma) \right).$$
\item[(iv) The point $\{\infty\}$] Change the domain of $ev$ by setting
$${\cal U}^\lambda =  \bigcup_{J \in \Jj} 
\mu^{\lambda}(A-Y,J) \times \mu(Y,J) \times_{G_{\infty} }S^2 $$
where $G_{\infty}$ is the four dimensional set of holomorphic maps from $S^2$ to itself which
fix $\{\infty\}$, and let
$$ev(\lambda, J, u,v,w_0) = \left( v(\infty), u(0), u(\infty), v(w_0),  
\int_D u^*(\omega \oplus \sigma)  \right).$$
\end{description}

The proofs of the transversality and dimension results for these cases  are analogous to the
case examined in Proposition \ref{proposition: DJempty}.   Hence, for each case there
is a set of second category of regular $J$ for which there will be no two component
bubbles of a certain homology class $Y$ which are not multiply covered.   Now, for each case,
take the $J$ which are in the set for all $Y$, i.e. the intersection over the countable set of possible 
homology classes $Y$.    The interesection of these five sets (one for each case) is
the set of regular almost complex structures described in the proposition.
  
To show that multiple bubbles would not occur, the argument from the
proof of Proposition \ref{proposition: DJempty} can be 
modified.  For each additional bubble, we would increase the number of homology classes 
used to form $\Uu_J$ by 1 and increase the number of $S^2$ used
by 2.  (See Theorem 6.3.2 from \cite{MS}). This adds $2n+2+4 = 2n+6$ to
the dimension of $\Uu_J$, and we may reduce by the six dimensional 
reparametrization group $G $ to get $2n$ added
dimensions.  The transversality results would carry through.  The codimension 
of $\Dd_J$ with one added bubble would increase by $2n+2$.  Hence,  
the codimension of $\Dd_J$ would be greater than the dimension of $\Uu_J$, 
so $\Dd_J$ will be empty.  
 
Finally, we must deal with the possibility of multiply covered curves. 
Without loss of generality, assume that the cusp curve has two
components: the $\lambda k(u)$ perturbed $J-$holomorphic component of class $A-Y$ 
and the $J-$holomorphic bubble component of class $X$.  Suppose
that $X$  has been reduced from the multiply covered $dX$ where 
$dX - Y = 0$ in homology for some positive integer $d > 1$.  Since
$M$ has dimension two or four, $M \times S^2$ has dimension less than or
equal to six.  Therefore, all
classes representable by a $J-$holomorphic or perturbed $J-$holomorphic 
curve give a nonnegative
integer when paired with the first Chern class.   In particular, $dX$ is representable so 
$$c_1(X) = \frac{1}{d} \cdot c_1(dX) \geq 0.$$  
This gives us $$c_1(A) = c_1(A-Y) + d \cdot c_1(X) > c_1(A-Y) + c_1(X) = c_1(A-Y+X).$$  
Consider the case fully explained in Proposition \ref{proposition: DJempty}; 
the others are identical.  When we imitate the proof of Lemma
\ref{lemma: dimUJ}, we see that the space we would consider as the domain
of the evaluation map is $\Uu_J' = \{ \lambda, \Uu_J^{\lambda '} \}$ where
$${\Uu}_J^{\lambda '} =  \mu^{\lambda}(A-Y,J) \times \mbox{ upper hemisphere of } S^2 \times 
\mu(X, J) \times _G S^2.$$
We calculate 
$$
\begin{array} {rcl}
\dim {\Uu}_J' & = &  1+(2c_1(A-Y) +2n +2) + 2  + (2c_1(X)+2n + 2) +2-6 \\
           & = & 3 + 2c_1(A-Y+X) + 4n\\
           & < & 3 + 2c_1(A) + 4n \\
           & = & 7 + 4n.
\end{array}$$
Note, however, that the codimension of  $\Dd_J$ will still be $7+4n$. Hence, again,
generically $\Dd_J$ will be empty.  We can deal with the case when $A-Y$ has
been reduced from a multiply covered component in a similar manner.
\QED

\section{Proofs of main theorems }
We will restate the theorems here as we prove them.

\begin{makingchzprod}
 Let $(M, \omega)$  be a compact symplectic manifold of dimension two or four. Then,
$$\cHZ (M \times D(a), \omega \oplus \sigma) \leq  a.$$
\end{makingchzprod}

\proof{}
Fix an almost complex structure $J \in \Jj$ on $M \times S^2$ so
that Proposition \ref{prop: nobub} holds.
Note that Proposition  \ref{prop: nobub} implies that 
Theorem \ref{theorem: HV1.12} and hence Theorem \ref{theorem: caHV}  
hold for $M$, if $M$ has dimension two or four.  This completes the proof.
\QED

\begin{makingmainthmgeod}
Let $(M, \omega)$ be a compact symplectic manifold  of dimension two or four. Let $\phi^H_t$ for
$0 \leq t \leq 1$ be a path in $Ham (M)$  generated by an autonomous 
Hamiltonian 
$H: M \rightarrow {\bf R}$  such that $\phi^H_0$ is the identity
diffeomorphism and  $\phi^H_t$ has no non-constant closed trajectory in
time less than 1.
Then, $\phi^H_t$ for $0 \leq t \leq 1$ is length minimizing among all 
homotopic paths between the identity and $\phi^H_1$.
\end{makingmainthmgeod}

\proof{} 
Theorem \ref{theorem: chzprod} implies that the capacity-area inequality
holds for $c_{HZ}$ for all split quasi-cylinders.  We can repeat the proof from
 Proposition 4.4 of \cite{LM2} to show that it holds for all quasi-cylinders.
Thus, ${\cHZ}$ satisfies condition (ii) of Theorem \ref{prop: criteria}  
for any Hamiltonian $H$ on $M$ if
$M$ has dimension two or four.  Now, we choose an autonomous  $H$ that generates
a flow $\phi_t^H$ which has no non-constant closed trajectories for $0 < t \leq 1$.  
In order to show that $H$ generates a path which is length minimizing  among
all homotopic paths, we must show that $c_{HZ}(H) \geq L(H)$ verifying
condition (i) of 
Theorem \ref{prop: criteria}.   
  We now invoke 
Proposition 3.1 from \cite{LM2}:

\begin{prop}
(Lalonde-McDuff) Let $M$ be any symplectic manifold and $H: M \rightarrow {\bf R}$ be
any compactly supported Hamiltonian with no non-constant closed 
trajectory in time less than 1.  Then
$${\cHZ}(H) \geq L(H).$$
\end{prop}
\proof{}
We give here a sketch of the proof.
Using $H$, we can construct a specific Hamiltonian $\overline{H}$ on $R_H^-(\frac{\nu}{2})$
and show that $\overline{H} \in {\cal H}_{ad}(R_H^-(\frac{\nu}{2})).$  Then, it
is easy to show that $m(\overline{H}) \geq m(H) = L(H)$, so 
 ${\cHZ}(R_H^-(\frac {\nu}{2})) \geq L(H)$ and hence ${\cHZ}(H) \geq L(H)$.
\QED

It follows that $\phi_t^H$ for $0 \leq t \leq 1$
is  length minimizing   among all paths homotopic to it with fixed endpoints
from the identity to $\phi_1^H$, and we are finished with the proof
of Theorem \ref{theorem: mainthmgeod}. \QED

\appendix
\section{Appendix}
Here we use Proposition \ref{prop: criteria} and $c_G$ to show two natural paths
in $Ham(\CP)$ and $Ham(\CPb)$ given by rotation are length minimizing.  The ball
embeddings are described explicitly.  Of course, Theorem \ref{theorem: rotCP2}
and Theorem \ref{theorem: rotCPtwoB} are also special cases of Theorem 
\ref{theorem: mainthmgeod}.

\subsection{Rotation in $\CP$}

\begin{theorem}
\label{theorem: rotCP2}
The path $\phi^P_t$ for $0 \leq t \leq 1$ in $Ham(\CP)$ given by
$$\phi^P _{t} [z_0: z_1: z_2] = [e^{\pi i t}z_0: z_1: z_2]$$
is length minimizing between the identity ($\phi^P_0$) and rotation
by $\pi$ radians in the first coordinate ($\phi^P_1$).  
\end{theorem}

\proof{}
To prove this theorem, we will use Gromov capacity $c_G$ and the criteria
from Proposition \ref{prop: criteria}.   To use these criteria, we need the
generating Hamiltonian of the path and its length.
The Hamiltonian function $P: \CP \rightarrow {\bf R}$ which generates our path
$\phi^P_t$ is 
$$P([z_0: z_1: z_2]) = \frac{\pi}{2} \frac{|z_0|^2}{|z_0|^2 + |z_1|^2
+ |z_2|^2}. $$ 

\begin{lemma}
\label{lemma: lengthofP}
The Hamiltonian $P$ defined on $\CP$ has $L(P) = \frac{\pi}{2}.$
\end{lemma}

\proof{}
Since $P$ is autonomous, 
$$L(P) = \max_{x \in \CP} P(x) - \min_{x \in \CP} P(x) = \frac{\pi}{2} - 0 = \frac{\pi}{2}.$$
\QED

Note that the criteria from Proposition \ref{prop: criteria} only
tells us if $\phi^P_t$ will be length minimizing within its homotopy class.
However, we use a proposition from \cite{LM2} to show it is actually
globally length minimizing.

\begin{proposition}
\label{proposition: r1}
(Lalonde-McDuff) Suppose we have a manifold $M$ and a capacity $c$ which satisfies 
condition (ii) of Proposition \ref{prop: criteria}.  The path
$\phi_t^H$ is length minimizing amongst all paths with the same endpoints
if $c(H) = L(H) \leq \frac{r_1(M)}{2}$.
\end{proposition}

The function $r_1$ is defined as follows:
if we let ${\cal L} : \pi_1(Ham^c(M)) \rightarrow {\bf R}$  be defined as
$${\cal L} ([\ga])  = \inf_{\ga \in [\ga]}L(\ga)$$
then
$$r_1(M) = \inf ( \{ {\rm Im} \; {\cal L} : \pi_1(Ham^c(M) \rightarrow {\bf R} \} \cap (0, \infty) )$$
if this set is not empty, and $\infty$ otherwise.

 Now, $\pi_1(Ham({\CP})) = {\bf Z}_3$,
generated by rotation through $2\pi$ radians in one coordinate \cite{P1},
specifically  the loop
$$\psi_t [z_0: z_1: z_2] = [e^{2\pi it}z_0: z_1: z_2]$$
for $0 \leq t \leq 1$.  By Theorem \ref{theorem: mainthmgeod}, 
the path $\psi_t$ for $0 \leq t \leq 1 - \epsilon$ is length minimizing
among homotopic paths for any $\epsilon > 0$.   By a limiting argument,
it is easy to see that $\psi_t$ for $0 \leq t \leq 1$ is also length minimizing
among homotopic paths.   Therefore, 
$${r_1}(\CP) = L(\psi_{t \in [0.1]}) = \pi.$$   
By Lemma \ref{lemma: lengthofP},
$L(\phi_t^P) = L(P) = \frac{\pi}{2} = \frac{r_1(M)}{2}$.   Hence, Proposition \ref{proposition: r1}
implies that if the hypotheses
from Proposition \ref{prop: criteria} are satisfied, $\phi^P_t$  
will actually be length minimizing among non-homotopic paths as well
as homotopic ones.

Since $L(P) = \frac{\pi}{2}$, we need to show $c_G(P) = \frac{\pi}{2}$.  Recall that the 
capacity of $P$ is the minimum of the capacities of $R_P^-$ and
$R_P^+$, the regions below and above the graph 
$$\Gamma_P = \{ x,s,t \; \vline \; P(x) = s \}$$
of $P$ in the six dimensional compact manifold
$$ {\CP} \times [0,\frac{\pi}{2}] \times [0,1]$$ 
endowed with the product symplectic form $\tau_0 \oplus ds \wedge dt$.

Since $R_P^-$ and $R_P^+$ are quasi cylinders with area $L(P) = \frac{\pi}{2}$
and the capacity area inequality holds on $\CP$ for $c_G$, we know
that 
$$c_G(P) \leq L(P) = \frac{\pi}{2}.$$  To show that $c_G(P) \geq \frac{\pi}{2}$,
we show both $c_G(R_P^-) \geq \frac{\pi}{2}$ and  $c_G(R_P^+) \geq \frac
{\pi}{2}$ by symplectically embedding a 6-ball of radius $1/ \sqrt{2} - \epsilon$ into
each region.  

We explicitly construct a symplectic embedding of $B^6(\frac{1}{\sqrt{2}}-\epsilon)$ into $R_P^-$ and $R_P^+$.
First, we consider 
$$R_P^- = \left \{ [z_0: z_1: z_2], s, t \; \vline \; 0 \leq t \leq 1,
\ell(t) \leq s \leq \frac{\pi}{2} \frac{|z_0|^2}{|z_0|^2 + |z_1|^2
+ |z_2|^2} \right \}$$
$$ \subset \CP \times [0, \frac{\pi}{2}] \times [0,1]$$
where $\ell(t)$ is some negative function close to 0.  In fact, we will
embed \linebreak
$B^6(\frac{1}{\sqrt{2}}-\epsilon)$ in the subset of $R_P^-$ where 
$s \geq 0$.    The embedding will be done in two steps: first, we embed
a 4-ball in $\CP$ and then embed a 2-ball in the two extra graph dimensions.

\begin{figure}[htb]
\centering
\begin{picture}(150,200)(40, -5)
\put(0,0){\includegraphics{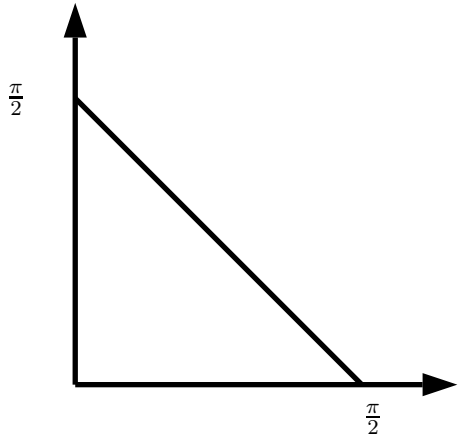}}
\put(25,125){$ \frac{\pi}{2}$}
\put(160,4){$\frac{\pi}{2}$}
\end{picture}
\caption{Image of $\CP$ under $\rho$ }
\label{figCP}
\end{figure}

To understand what is happening geometrically, we  identify
$\CP$ with its image under the moment map of the 
$T^2$ action 
$$(\theta_0, \theta_1) ([z_0: z_1: z_2]) = [e^{\pi i \theta_0}z_0:
e^{ \pi i \theta_1}z_1: z_2]$$
with $0 \leq \theta_0, \theta_1 \leq 1$.  The moment map for this action ${\rho} : {\CP} \rightarrow {\bf R^2}$ is given by
$${\rho} ([z_0: z_1: z_2] ) =
 \left( \frac{\pi}{2} 
\frac{ \vline z_0 \vline ^2}{\vline z_0 \vline ^2 + \vline z_1 \vline ^2
+ \vline z_2 \vline ^2},  \frac{\pi}{2} 
\frac{ \vline z_1 \vline ^2}{\vline z_0 \vline ^2 + \vline z_1 \vline ^2
+ \vline z_2 \vline ^2} \right) $$
and the image of $\CP$ under ${\rho}$ is the right triangle pictured
in Figure \ref{figCP}.  The Hamiltonian $P$ is projection onto the horizontal axis 
and its image is the interval $[0, \frac{\pi}{2}]$.

\begin{figure}[htb]
\centering
\begin{picture}(150,200)(40, -5)
\put(0,0){\includegraphics{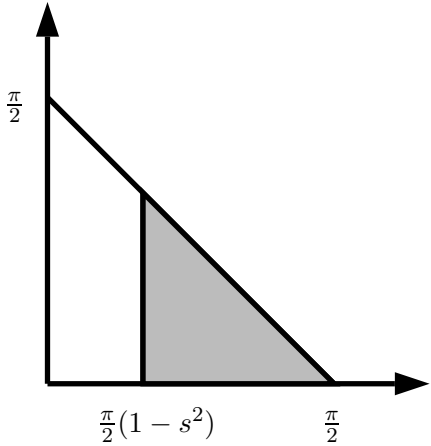}}
\put(40,125){$\frac{\pi}{2}$}
\put(160,4){$\frac{\pi}{2}$}
\put(75,4){$\frac{\pi}{2}(1-s^2)$}
\end{picture}
\caption{Image of $B^4(s)$ under ${\rho} \circ i^-$}
\label{figimageofi-}
\end{figure}

Let $i^-: {\bf C^2} \rightarrow \CP$ be the map 
$$i^-(z_1, z_2) = [\sqrt{1 - |z_1|^2 - |z_2|^2}: z_1: z_2].$$
Note that $i^-$ restricted to $B^4(s) = \{ (z_1, z_2) \; \vline \; |z_1|^2 +
|z_2|^2 \leq s^2\}$ is a symplectic embedding for $s < 1$.  The image
of $i^-$ composed with ${\rho}$ is the shaded triangle in Figure
\ref{figimageofi-}.

\begin{figure}[htb]
\centering
\begin{picture}(300,300)(90, -15)
\put(0,0){\includegraphics{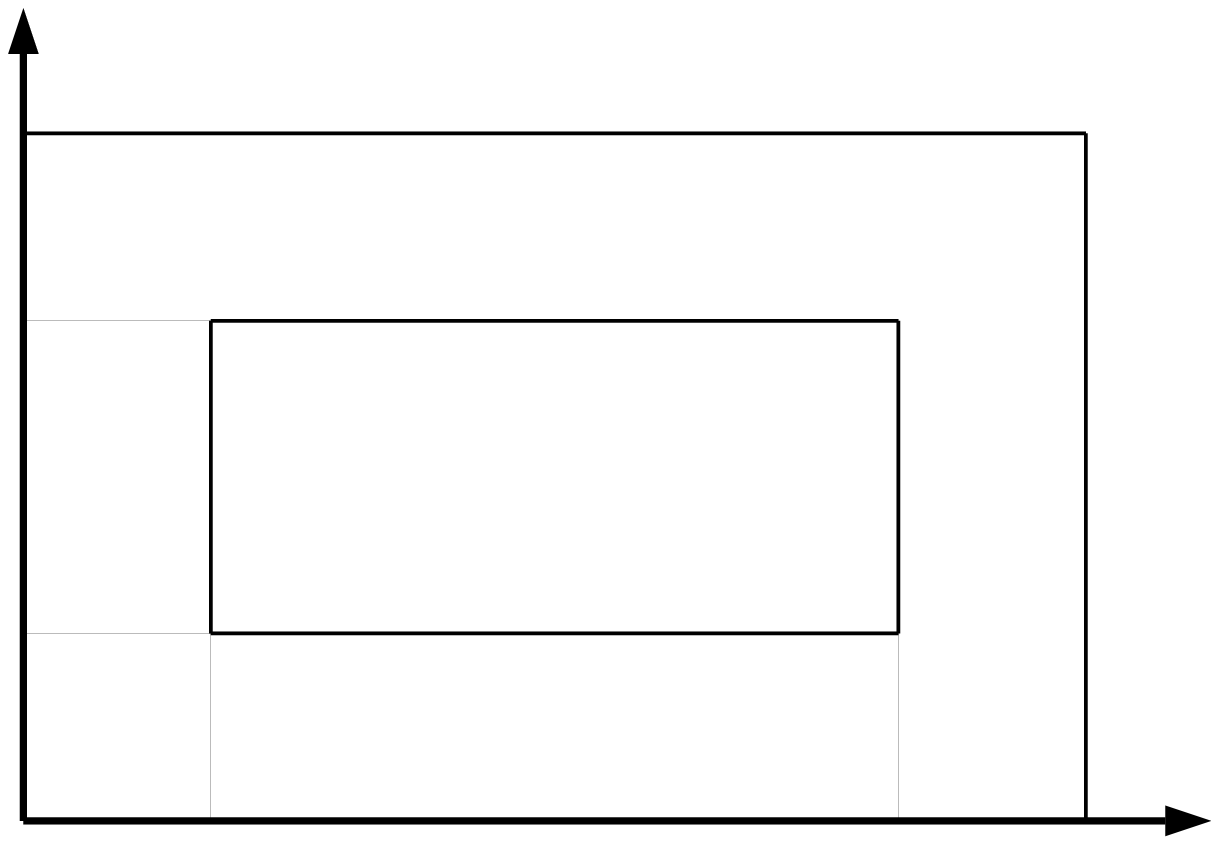}}
\put(60, 220){$1$}
\put(45, 165){$\frac{1}{2} + \frac{1}{\sqrt{2}} r$}
\put(45,  75){$\frac{1}{2} - \frac{1}{\sqrt{2}} r$}
\put(60, 20){$0$}
\put(105, 5){$\frac{\pi}{4} + (\frac{\pi}{2})r^2 - (\frac{\pi}{\sqrt{2}})r$}
\put(320, 5){$\frac{\pi}{4} + (\frac{\pi}{2})r^2 $}
\put(395,5){$\frac{\pi}{2}$}
\end{picture}
\caption{Image of $B^2(r)$ under $\psi_r^-$}
\label{figCP2rect-}
\end{figure}

Choose an $r \leq 1 / \sqrt{2}$.  For any $\epsilon > 0$,  we can symplectically embed $B^2(r - \epsilon)$ (the closed
2-ball of radius $r - \epsilon$) into the smaller rectangle in Figure \ref{figCP2rect-}
because the area of the ball is $\pi (r- \epsilon)^2$ and the area of the rectangle is
$(\frac{\pi}{\sqrt{2}}r)(\frac{2}{\sqrt{2}}r) = \pi r^2$.  Denote this
mapping by $\psi_r^-$.  Let $R = \frac{1}{\sqrt{2}} - \epsilon$.  It is possible
to choose the family of maps $\psi_r^-$ so that they fit together
to form a smooth map $\psi_R^-$ on $B^2(R)$ such that for $r < R$, 
$$\psi_R^- \vline_{B^2(r)} = \psi_r^-.$$
In particular,  this means the images of nested circles under $\psi_r^-$
are disjoint and nested inside the larger 
rectangle in Figure \ref{figCP2rect-}.  

We define the map $\Psi^-: B^6(\frac{1}{\sqrt{2}}- \epsilon) \rightarrow R_P^-$ by 
$$ \Psi^-(z_1, z_2, u, v) = (i^-(z_1, z_2), \psi_{R}^-(u,v))$$
where the domain coordinates lie in ${\bf C} \times {\bf C} \times {\bf R}
\times {\bf R}$ and satisfy $|z_1|^2 + |z_2|^2 + u^2 + v^2 \leq (1/\sqrt{2} - \epsilon)^2$.
$\Psi^-$ will be the required embedding.  We must show that $\Psi^-$ is well
defined, i.e. the image of $\Psi^-$ does actually lie in $R_P^-$.  Once this
has been demonstrated, it is easy to see that $\Psi^-$ is symplectic since it
is the product of two symplectic maps into a symplectic manifold given the
product symplectic structure.   

Since the map $i^-$ obviously is a well defined embedding, we must only check that
for a given point $(z_1, z_2, u, v) \in B^6(\frac{1}{\sqrt{2}} - \epsilon)$,
the image of $\psi_R^- (u,v)$ is contained in
$[0, \frac{\pi}{2}(1 - |z_1|^2 - |z_2|^2)] \times [0,1] \subset {\bf R^2}$.
We let $u^2 + v^2 = r^2$ and use the fact that
$$\psi_R^- \vline_{B^2(r)} = \psi_r^-.$$  The height of the rectangle
(the second coordinate of the image of $\psi_r^-$) covers the region 
$$\left[ \frac{1}{2} -
\frac{1}{\sqrt{2}}r, \frac{1}{2} + \frac{1}{\sqrt{2}}r \right]$$ 
which is contained in the required interval $[0,1]$ for
all $r \in [0, \frac{1}{\sqrt{2}}]$.  For any given $r$, the width of the
rectangle (the first coordinate of the image of $\psi_r^-$) covers the
region 
$$\left[ \frac{\pi}{4} + \frac{\pi}{2}r^2 - \frac{\pi}{\sqrt{2}}r,
\frac{\pi}{4} + \frac{\pi}{2}r^2 \right].$$
As is required, the function $ \frac{\pi}{4} + \frac{\pi}{2}r^2 - \frac{\pi}{\sqrt{2}}r$ is
greater than zero and decreasing for all values of $r \in [0, \frac{1}{\sqrt{2}}]$.  For the final check, we must examine the upper endpoint of the
first coordinate of the image of $\psi_r^-$, $\frac{\pi}{4} + \frac{\pi}{2}r^2$,
to ascertain that it is less than
or equal to $\frac{\pi}{2}(1 - |z_1|^2 - |z_2|^2) = P \circ i^- (z_1, z_2)$ for $u,v$ such that
$(z_1, z_2, u, v) \in B^6(\frac{1}{\sqrt{2}} - \epsilon)$.
This is a simple calculation hinged on the fact that $P$ applied to the 
image under $i^-$ of any 3-sphere is constant:
$$\begin{array} {rcl}
\frac{\pi}{4} + \frac{\pi}{2}r^2 & = &\frac{\pi}{4} + \frac{\pi}{2}(u^2 + v^2) \\
                                 & \leq & \frac{\pi}{4} + \frac{\pi}{2}
(\frac{1}{2} - |z_1|^2 - |z_2|^2) \\
                                 & = & \frac{\pi}{2}(1 - |z_1|^2 -|z_2|^2).
\end{array}$$
Hence, the map $\Psi^-$ is a well defined symplectic embedding of $B^6(\frac
{1}{\sqrt{2}} - \epsilon)$ into $R_P^-$.

In a similar manner we can define an embedding $\Psi ^+$ of 
$B^6(\frac{1}{\sqrt{2}} - \epsilon)$
into $R_P^+$, the region above $\Gamma_P$.  We do this in two parts,
as before, but now we want to center our ball in the $\CP$ portion away
from $[1:0:0]$.  

\begin{figure}[h]
\centering
\begin{picture}(150,200)(40, -5)
\put(0,0){\includegraphics{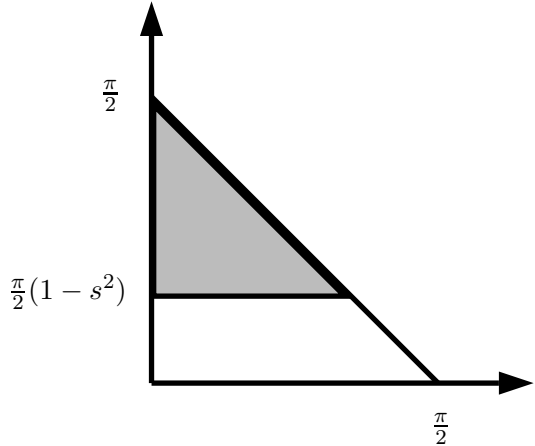}}
\put(60,130){$\frac{\pi}{2}$}
\put(185,4){$\frac{\pi}{2}$}
\put(25,55){$\frac{\pi}{2}(1-s^2)$}
\end{picture}
\caption{Image of $B^4(s)$ under ${\rho} \circ i^+$}
\label{figimageofi+}
\end{figure}

Let $i^+: {\bf C^2} \rightarrow \CP$ be the map 
$$i^+(z_1, z_2) = [z_1: \sqrt{1 - |z_1|^2 - |z_2|^2}: z_2].$$
Note that $i^+$ restricted to $B^4(s) = \{ z_1, z_2 \; \vline \; |z_1|^2 +
|z_2|^2 \leq s^2\}$ is a symplectic embedding for $s < 1$.
The image of $i^+$ composed with ${\rho}$ is the shaded triangle in Figure \ref{figimageofi+}.

The map $i^+$ will be the first part of $\Psi^+$.

\begin{figure}[htb]
\centering
\begin{picture}(300,300)(90, -15)
\put(0,0){\includegraphics{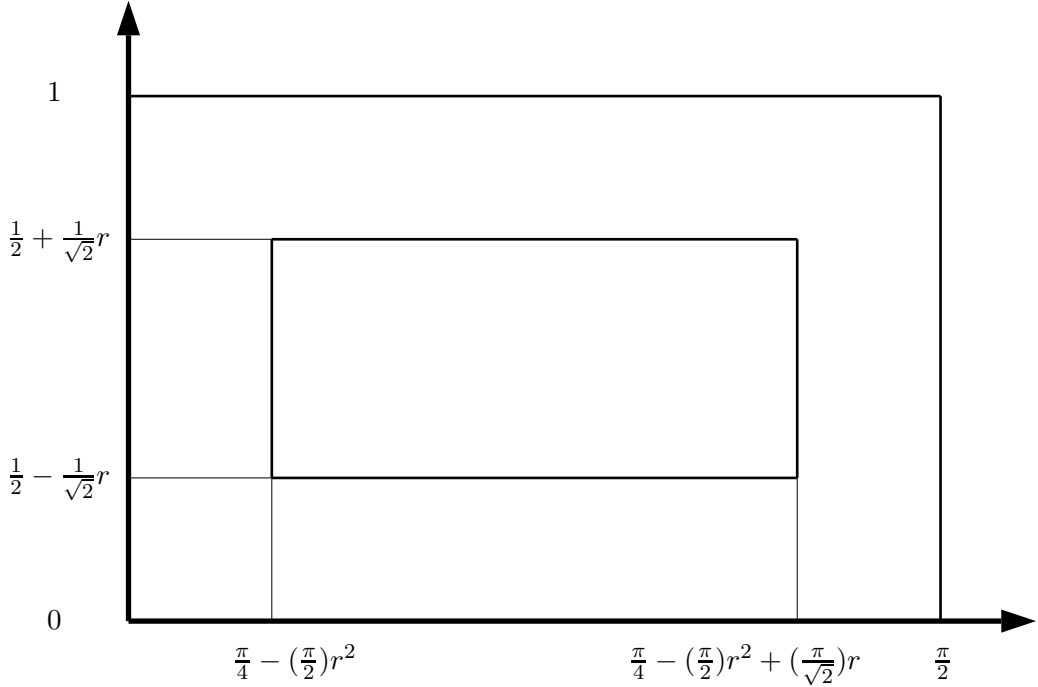}}
\put(60, 220){$1$}
\put(45, 165){$\frac{1}{2} + \frac{1}{\sqrt{2}} r$}
\put(45,  75){$\frac{1}{2} - \frac{1}{\sqrt{2}} r$}
\put(60, 20){$0$}
\put(280, 5){$\frac{\pi}{4} - (\frac{\pi}{2})r^2 + (\frac{\pi}{\sqrt{2}})r$}
\put(130, 5){$\frac{\pi}{4} - (\frac{\pi}{2})r^2 $}
\put(395,5){$\frac{\pi}{2}$}
\end{picture} 
\caption{Image of $B^2(r)$ under $\psi_r^+$}
\label{figCP2rect+}
\end{figure}

Next, note that we can symplectically embed $B^2(r - \epsilon)$ into the  smaller
rectangle in Figure \ref{figCP2rect+} because the area of the ball is 
$\pi (r-\epsilon)^2$ and the area of this rectangle is $\pi r^2$. 

 We denote this mapping by $\psi_r ^+$.
As in the previous set up, we may assume that for $r < R$, $\psi_R ^+ \vline_{B^2(r)} = \psi_r ^+$.  Then, we define
$\Psi ^+: B^6(\frac{1}{\sqrt{2}} - \epsilon) \rightarrow R_P^+$ by 
$$ \Psi ^+ (z_1, z_2, u, v) = (i^+ (z_1, z_2), \psi_R ^+ (u,v))$$
where the domain coordinates lie in ${\bf C} \times {\bf C} \times {\bf R}
\times {\bf R}$ and satisfy $|z_1|^2 + |z_2|^2 + u^2 + v^2 \leq (1/\sqrt{2} - \epsilon)^2$.
Just as we checked that $\Psi$ is a well defined symplectic embedding,
we may verify that $\Psi^+$ is also a well defined symplectic
embedding. \QED

\subsection{Rotation in $\CPb$}
The next natural path to examine is rotation on the symplectic
blow-up $\CPb$ of $\CP$.   For precise details of the definition of $\CPb$, 
see Chapter Six of \cite{MSintro}.  Geometrically,
$\CPb$ can be thought of as  
 the manifold obtained by removing from $\CP$ an open 4-ball of radius
$\lambda$ centered at $[1:0:0]$ and collapsing its boundary $S^3$ along the fibers 
of the Hopf map.  The collapsed $S^3$, now an $S^2$, is the
exceptional divisor $\Sigma$:
$$\Sigma = \left \{ [z_0:z_1:z_2] \; \vline \; \frac{\vline z_1 \vline ^2 + 
\vline z_2\vline ^2}{\vline z_0 \vline ^2} = \lambda ^2 \right \} / \sim$$
where 
$$[z_0:z_1:z_2] \sim [w_0:w_1:w_2] \mbox{ if } z_1 w_2 = z_2 w_1. $$
 This is the interpretation of $(\CPb, \tau_\lambda)$ most often referred to in this paper.
   
Alternatively, if we think of $\CP$ as a 4-ball of radius 1 with the boundary $S^3$
collapsed along the fibers of the Hopf map, then $(\CPb, \tau_{\lambda})$ is an
annulus 
$\{ (w_0, w_1) \; \vline \; (1-\lambda^2) \leq |w_0|^2 + |w_1|^2 \leq 1\}$
with both boundaries collapsed along the Hopf fibers.  

The rotation $\phi^P _t$ on $\CP$ in the first homogeneous coordinate descends to a
well defined rotation on $\CPb$.  To check this, it is only necessary to verify
that the rotation keeps the set of removed points invariant and that the 
rotation is  well defined under the equivalence imposed on the boundary. 
In the same way, one can see that the projection of rotation in the second 
homogeneous coordinate in $\CP$
$$\phi^Q _t [z_0: z_1: z_2] = [z_0: e^{\pi it}z_1: z_2]$$ is well defined
on $\CPb$.    It is important to realize, however,  that the rotation
in the first homogeneous coordinate is qualitatively different from rotation
in the second:
$\phi ^P _t$ fixes each point on $\Sigma$, whereas $\phi ^Q _t$ 
keeps $\Sigma$ invariant but the points on $\Sigma$ rotate.

Note that the function $P$ is well defined on $\CPb$.  When blowing up, we 
collapse the boundary of the ball of radius $\la$ along orbits of an $S^1$ 
action, and $P$ (defined on $\CP$) is invariant under this action. 
Therefore, $P$, defined appropriately, is the Hamiltonian function which 
generates rotation in the first homogeneous coordinate on 
$\CP$ and $\CPb$.  Hence,
we will use $P$ to denote this Hamiltonian function and $\phi ^P _t$ to
denote its flow on both manifolds,  and
it will be clear from context which domain we are considering.  Similarly,
if we let $Q: \CP \rightarrow {\bf R}$ be defined as
$$ Q[z_0: z_1: z_2] = \frac{\pi}{2} 
\frac{ \vline z_1 \vline ^2}{\vline z_0 \vline ^2 + \vline z_1 \vline ^2
+ \vline z_2 \vline ^2}$$
it is clear that $Q$ induces the rotation in the second coordinate 
$\phi ^Q _t$ on $\CP$ and $\CPb$.  

\subsubsection{Rotation in $\CPb$ induced by $P$}

Here we begin our treatment of the rotation induced by $P$ applied to $\CPb$. 
In moving from $\CP$ to $\CPb$,
we have altered the domain of $P$ in a consequential way.

\begin{lemma}
\label{lemma: lengthofPCPb}
The Hamiltonian $P$ defined on $\CPb$ has $L(P) = \frac{\pi}{2}(1 - \la^2)$.
\end{lemma}
\proof{}
Written out in homogeneous coordinates, 
$$\CPb = \{ [ \sqrt{1- |z_1|^2 - |z_2|^2}: z_1: z_2] \; \vline \; 
\lambda^2 \leq |z_1|^2 + |z_2|^2 \leq 1\}$$
with the appropriate equivalence relation on the exceptional divisor.  
Hence, it is easy to see that 
$$L(P) = \max_{x \in \CPb} P(x) - \min_{x \in \CPb} P(x) = \frac
{\pi}{2}(1 - \lambda^2) - 0 = \frac{\pi}{2}(1 - \lambda^2) .$$
\QED

\begin{figure}[htb]
\centering
\begin{picture}(150,200)(40, -5)
\put(0,0){\includegraphics{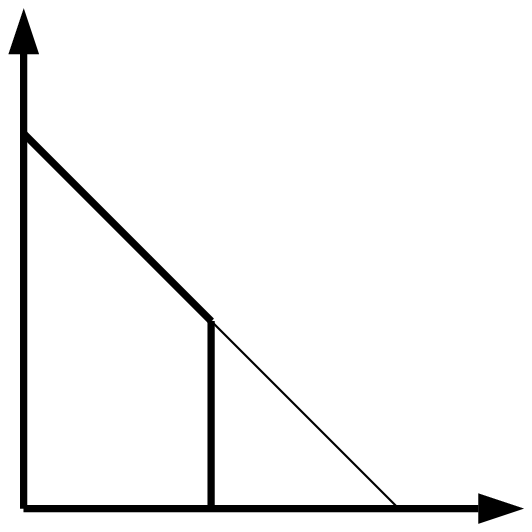}}
\put(35,125){$\frac{\pi}{2}$}
\put(160,4){$\frac{\pi}{2}$}
\put(80,4){$\frac{\pi}{2}(1-\lambda^2)$}
\end{picture} 
\caption{Image of $\CPb$ under ${\rho}$ }
\label{figCPb0}
\end{figure}

The image of $\CPb$ under the map ${\rho}$ is the quadrilateral
in Figure \ref{figCPb0}.   
Since the map $P$ is projection onto the
horizontal axis, with this quadrilateral as its domain, $P$
has image $[0, \frac{\pi}{2}(1-\la^2)]$.  This verifies Lemma
\ref{lemma: lengthofPCPb}.

\begin{theorem}
\label{theorem: rotCPtwoB}
The path $\phi^P_t$ for $0 \leq t \leq 1$ in $Ham(\CPb)$ given by
$$\phi^P _{t} [z_0: z_1: z_2] = [e^{\pi i t}z_0: z_1: z_2]$$
is length minimizing between the identity ($\phi^P_0$) and rotation
by $\pi$ radians in the first coordinate ($\phi^P_1$).  
\end{theorem}

\proof{}

By using the embeddings from the $\CP$ case adjusted appropriately, we can show  that
$c_G(P) = \frac{\pi}{2} (1 - \la^2)$.  By Proposition \ref{prop: criteria}
and Lemma \ref{lemma: lengthofPCPb},
this will tell us that $\phi^P_t$ is length minimizing in its homotopy
class.   Proposition \ref{proposition: r1} can be applied in the same way
as in the proof of Theorem \ref{theorem: rotCP2} to show that $\phi^P_t$
is actually globally length minimizing.   We omit the details but note
that  $\pi _1(Ham(\CPb)) = {\bf Z}$  generated by the loop $\psi_t$
described in the first section of this paper, see
\cite{AM} and \cite{P2}. 

To show that  
$c_G(R_P^+) \geq \frac{\pi}{2}(1 - \lambda^2)$ requires no additional
work; we may use the embedding $\Psi^+$ from the $\CP$ case. 
However, to prove \linebreak
$c_G(R_P^-)  \geq \frac{\pi}{2}(1 - \lambda^2)$ takes some manipulation.
We must produce a new embedding 
$\Upsilon^- : B^6(\sqrt{ \frac{1-\lambda^2}{2}}
- \epsilon) \rightarrow R_P^-$ because the old embedding, $\Psi^-$, has in
its image some points that were removed under the blow-up.

\begin{figure}[htb]
\centering
 \begin{picture}(150,300)(40, -5)
\put(0,0){\includegraphics{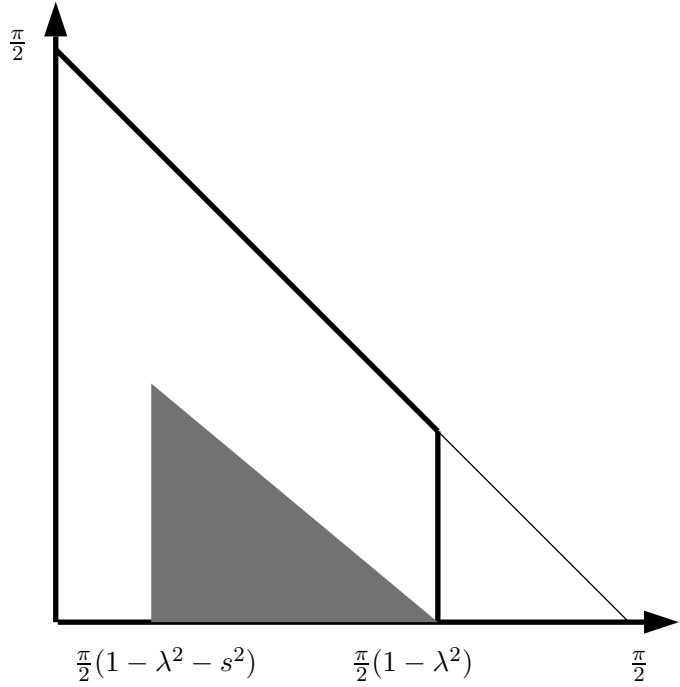}}
\put(20,260){$\frac{\pi}{2}$}
\put(255,25){$\frac{\pi}{2}$}
\put(45,25){$\frac{\pi}{2}(1- \lambda^2 - s^2)$}
\put(150, 25){$\frac{\pi}{2}(1- \lambda^2)$}
\end{picture} 
\caption{Image of $B^4(s)$ under ${\rho} \circ j^-$}
\label{figimageofj-}
\end{figure}

Consider the open shaded triangle in Figure \ref{figimageofj-} for some
$s$ where 
$s^2 \in [0, 1-\lambda^2]$.

By Delzant's theorem, the preimage of this set under the map
${\rho}$  is a symplectic submanifold.  This preimage is
equal to the set $U_s \subset \CPb$ where
$$U_s = \left \{ [z_0: z_1: z_2] \; \vline \; |z_0|^2 = (1- \lambda^2 - \tau^2),
|z_1|^2 < \tau^2, 0 \leq \tau^2 < s^2 \right \}. $$

We will prove that there exists a symplectic embedding $j_s^-$ of 
$B^4(s - \epsilon)$ into $U_s$.
$U_s$ is symplectomorphic to
the set $V_s \subset {\bf R^4}$ where 
$$ V_s = \left \{ (z_0, B^2(\sqrt{1- \lambda^2 - |z_0|^2})) \; \vline \; 
z_0 \in {\bf C}, 1 - \lambda^2 - s^2 < |z_0|^2 < 1 - \lambda^2 \right \} .$$
$V_s$ is just a set of 2-balls fibered over an annulus.  If we cut this 
annulus to make it a rectangle (this does not change the symplectic
capacity), we arrive at the set
$$T_s = \left \{ (x,y,B^2(\sqrt{s^2-y^2})) \; \vline \; 
0 \leq x < \pi s, 0 \leq y < s \right \} \subset {\bf R^4}.$$
$T_s$ is a generalized trapezoid, that is it consists of balls fibered over a rectangle.  It is not hard to show that the capacity of $T_s$ is the same
as the capacity of the more standard trapezoid
$$T^4(\pi s^2) = \left \{  (x,y,B^2(\sqrt{s^2 - \frac{y}{\pi}}) \; \vline
\; 0 \leq x < 1, 0 \leq y < \pi s^2 \right \}.$$
In Lemma 3.6 of \cite{LM2}, it is shown that the capacity of 
$T^4(\pi s^2)$ is equal to the capacity of $B^4(s)$.  Hence,
$$c_G(U_s) = c_G(T_s) = c_G(T^4(\pi s^2)) = c_G(B^4(s)) = \pi s^2$$
and we can embed
$B^4(s - \epsilon)$ into $U_s$ for any
$\epsilon > 0$. Call this embedding $j_s^-$.  Consider the family
of maps $j_s^- : B^4(s - \epsilon) \rightarrow U_s$ for all $0 \leq s \leq \sqrt{1-\la ^2} $.  Let 
 $$S = \sqrt{\frac{1 - \lambda^2}{2}} -\epsilon.$$   
Without loss of generality, we may assume that the family of maps
$j_s$ satisfies
$$j_S^- \vline_{B^4(s)} = j_s^- $$ 
for $s \leq S$, so that 3-spheres of constant radius appear as vertical lines
in the moment map picture.  To be precise, if $(w_o, w_1) \in {\bf C^2}$
and $|w_0|^2 + |w_1|^2 = s^2$,
then ${\rho} \circ j_S^- (w_0, w_1)$ lies on the vertical line through
the point $(\frac{\pi}{2}(1- \lambda^2 - s^2), 0)$.  Thus, $P$ applied
to the image of 3-spheres under $j_S^-$ is constant.  

Now, we have an embedding $j_S^-$ from 
$B^4(\sqrt{\frac{1-\lambda^2}{2}} - \epsilon)$ into $\CPb$.
Our next task is to work with the other two dimensions and construct
$\Upsilon^-$.

\begin{figure}[htb]
\centering
\begin{picture}(300,300)(90, -15)
\put(0,0){\includegraphics{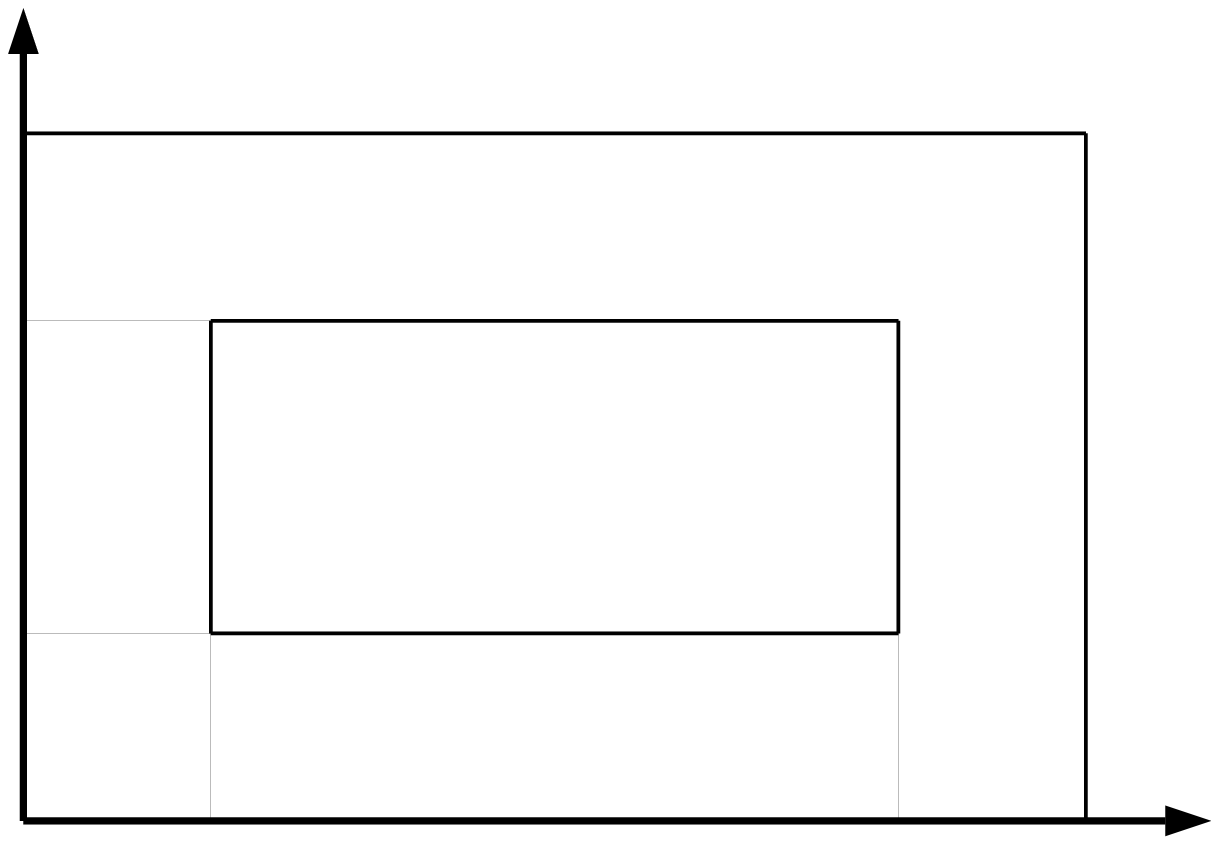}}
\put(60, 220){$1$}
\put(25, 165){$\frac{1}{2} + \frac{1}{\sqrt{2-2\lambda^2}} r$}
\put(25,  75){$\frac{1}{2} - \frac{1}{\sqrt{2-2\lambda^2}} r$}
\put(60, 20){$0$}
\put(95, 5){$\frac{\pi}{4}(1-\lambda^2) + \frac{\pi}{2}(r^2) - \frac{\pi}{\sqrt{2}}(1-\lambda^2)r$}
\put(285, 5){$\frac{\pi}{4}(1-\lambda^2) + \frac{\pi}{2}(r^2) $}
\put(395,5){$\frac{\pi}{2}(1-\lambda^2)$}
\end{picture} 
\caption{Image of $B^2(r)$ under $\upsilon_r^-$}
\label{figCP20rect-}
\end{figure}

Fix an $r < \sqrt{\frac{1-\lambda^2}{2}}$.  We can symplectically embed
$B^2(r)$ into the smaller rectangle in Figure \ref{figCP20rect-}
because the area of the ball is $\pi r^2$ and the area of the rectangle
is 
$$\left( \frac{\pi}{2} (\sqrt{1- \lambda^2}) r \right) 
\left( \frac{2}{\sqrt{2(1-\lambda^2)}} \right)
= \pi r^2.$$

Denote this embedding by $\upsilon_r^-$.  As before, we assume that
for $r < S$, $\upsilon_{S} ^-\vline_{B^2(r)} = \upsilon_r^-$, and define
$\Upsilon^-: B^6(S) \rightarrow R_P^-$
by  
$$\Upsilon^- (w_0, w_1, u, v) = \left( j^-(w_0, w_1), \upsilon_{S}^-(u,v) \right).$$
Using the fact that $P$ is constant along the image under $j^-$ of
3-spheres, it is routine to check that in fact $\Upsilon^-$ is well defined. \QED

\subsubsection{Rotation in $\CPb$ induced by $Q$}

Now, recall the Hamiltonian function $Q: \CP \rightarrow {\bf R}$ given by
$$Q([z_0: z_1: z_2]) = \frac{\pi}{2} \frac {|z_1|^2}
   {|z_0|^2 + |z_1|^2 + |z_2|^2}.$$
It  is easy to check that $L(Q) = \frac{\pi}{2}$ and that the flow of
$Q$ is the path in $Ham (\CP)$
$$\phi^Q _{t} [z_0: z_1: z_2] = [z_0: e^{\pi i t} z_1: z_2].$$ 
$Q$ descends to a well defined function on $\CPb$, and its flow
descends to a well defined rotation on $\CPb$.  

\begin{lemma}
\label{lemma: lengthofQCPb}
The Hamiltonian $Q$ defined on $\CPb$ has $L(Q) = \frac{\pi}{2}$.
\end{lemma}
\proof{}
Written out in homogeneous coordinates, 
$$\CPb = \{ [ \sqrt{1- |z_1|^2 - |z_2|^2}: z_1: z_2] \; \vline \; 
\lambda^2 \leq |z_1|^2 + |z_2|^2 \leq 1\}$$
with the appropriate equivalence relation on the exceptional divisor.  
Hence, it is easy to see that 
$$L(Q) = \max_{x \in \CPb} Q(x) - \min_{x \in \CPb} Q(x) = \frac
{\pi}{2} - 0 = \frac{\pi}{2}.  \mbox{\QED} $$

Recall that the image of ${\rho}$ applied to $\CPb$ is
the quadrilateral depicted in Figure \ref{figCPb0}.  The map $Q$ defined
on $\CPb$ is projection onto the vertical axis in this picture and has image 
$[0, \frac{\pi}{2}]$, verifying Lemma \ref{lemma: lengthofQCPb}.

We could show that the path $\phi ^Q _t$  for $0 \leq t \leq 1$
defined on $\CP$ is length minimizing
by using the argument from Theorem \ref{theorem: rotCP2}.  However, we cannot use the
arguments from Theorem \ref{theorem: rotCPtwoB} to show that $\phi ^Q _t$ is length minimizing
on $\CPb$.  Lemma \ref{lemma: lengthofQCPb} tells us that the length of $Q$ 
does not decrease when going from $\CP$ to  $\CPb$.   However, the volume
of the manifold $\CPb$ is less than the volume of $\CP$.
There is not a straight forward way to embed large enough 6-balls to 
show that $c_G(Q) = \frac{\pi}{2}$
on $\CPb$, i.e. 6-balls of raduis $\frac{1}{\sqrt{2}}$. 
(Recall that in the proof of Theorem \ref{theorem: rotCPtwoB} 
for $\phi^P_t$ we only
had to embed balls of radius close to $\sqrt{\frac{1- \lambda ^2}{2}}$.)
 
These two rotations, $\phi ^P_t$ and $\phi ^Q _t$,\ are essentially 
the only two different types of rotations of $\CP$ which descend to
rotations on $\CPb$.   In order for any  rotation
to descend properly from $\CP$ to $\CPb$, $[1: 0: 0]$ must be a fixed 
point of the
rotation in $\CP$.  A rotation in $\CP$ has one isolated fixed point and
a fixed sphere, e.g. $\phi ^P _t$ has an isolated fixed point at $[1: 0: 0]$
and a fixed sphere consisting of the points of the form $[0: z_1: z_2] 
\subset \CP$.   Since $[1: 0: 0]$ is the isolated fixed point of
$\phi ^P _t$ and a point on the fixed sphere of $\phi ^Q _t$, we have
accounted for both types of rotations.  

Because $L(Q)$ does not decrease when moving from $\CP$ to $(\CPb, \tau_{\lambda})$,  
we cannot use the Gromov capacity to show that the
rotation induced by $Q$ on $\CPb$ is length minimizing.

\medskip
\medskip
\noindent \textsc{D\'epartement de Math\'ematiques, 
        Universit\'e du Qu\'ebec \`a Montr\'eal, \linebreak
        CP 8888, Succ. Centre Ville,
        Montr\'eal, Qu\'ebec   H3X 3P8,
        Canada}

\medskip
\noindent \textit{E-mail address:}  \texttt{jennifer@math.uqam.ca}

\end{document}